\documentclass{amsart}
\usepackage{amsmath,amsfonts,amssymb,amsthm,mathtools,enumitem}
\usepackage{geometry}
\usepackage{microtype}
\usepackage{hyperref}

\newtheorem{theorem}{Theorem}[section]
\newtheorem{lemma}[theorem]{Lemma}
\newtheorem{proposition}[theorem]{Proposition}
\newtheorem{corollary}[theorem]{Corollary}
\theoremstyle{definition}
\newtheorem{definition}[theorem]{Definition}
\newtheorem{remark}[theorem]{Remark}

\DeclarePairedDelimiter{\abs}{\lvert}{\rvert}
\DeclarePairedDelimiter{\norm}{\lVert}{\rVert}
\DeclareMathOperator{\divg}{div}
\newcommand{\vu}{\mathbf{u}}
\newcommand{\vv}{\mathbf{v}}
\newcommand{\R}{\mathbb{R}}
\newcommand{\E}{\mathbb{E}}
\newcommand{\Ray}{\mathrm{Ra}}
\newcommand{\Nuss}{\mathrm{Nu}}
\newcommand{\ez}{\mathbf{e}_z}
\newcommand{\LL}{\mathcal{L}}
\numberwithin{equation}{section}
\begin{document}

\title{A lower bound on the Nusselt number in stochastically forced Rayleigh--B\'enard convection}
\renewcommand{\shorttitle}{A lower bound on the Nusselt number}

\author{Juraj Földes}
\address{Department of Mathematics, University of Virginia, Kerchof Hall,
  Charlottesville, VA 22904-4137}
\email{foldes@virginia.edu}
\thanks{J. F. was partially supported by the grants NSF-DMS-1816408, Simons foundation SFI-MPS-TSM-00013942, and U.S. DOE DE-NA0004269
}

\author{Jun Hyun Park}
\address{Department of Mathematics, University of Virginia, Kerchof Hall,
  Charlottesville, VA 22904-4137}
\email{rsx5wp@virginia.edu}

\date{}

\begin{abstract}
We establish a lower bound on the expected long-time averaged Nusselt number
for two-dimensional Rayleigh--B\'enard convection at infinite
Prandtl number, with free-slip boundaries and a suitable additive noise of
amplitude $0<\sigma\leq1$ in the temperature equation. For every
deterministic initial temperature, we prove that
\[
    \underline{\Nuss}-1
    \geq c(\sigma^3\Ray)^{7/192},
\]
provided that
\[
    \sigma^3\Ray
    \geq C(1+\log(1/\sigma))^{192/7},
\]
where $c,C>0$ depend only on the domain parameter $L$.
Without noise, no such bound can hold, since the conductive
state is a solution with $\Nuss=1$ for every Rayleigh number.

The proof uses the auxiliary functional method. For general
Markov processes, we show that the defining inequality of an
auxiliary functional need only hold after a short time average.
This reduces the problem to showing that, with uniformly positive probability, the solution quickly leaves the set of states with low heat flux and low thermal dissipation. We establish
this by projecting onto an unstable eigenvector of the
linearization around the initial mean temperature profile,
yielding a one-dimensional It\^o process with linear growth,
additive noise, and nonlinear error terms.



\end{abstract}
\maketitle
\tableofcontents
\section{Introduction}\label{sec:intro}\label{subsec:main-result}

Rayleigh--B\'enard convection describes the flow of a fluid heated from below and cooled from above. It arises in a variety of physical settings, including geophysics, oceanography, and atmospheric science. At sufficiently large thermal forcing, the resulting flow can become turbulent, exhibiting complex spatial structures and chaotic temporal dynamics. Rather than attempting to characterize individual solutions, one is often interested in their long-time statistical properties. A central observable in this context is the Nusselt number, $\Nuss$, which measures the long-time average of vertical heat transport relative to the heat transported by pure conduction. Thus, $\Nuss=1$ corresponds to a fluid at rest, while $\Nuss>1$ indicates enhanced heat transport due to convection.

A fundamental question in the theory of thermal convection is how $\Nuss$ scales with the Rayleigh number $\Ray$, which measures the strength of buoyancy relative to viscous damping and thermal diffusion. In many cases of interest, such as atmospheric or ocean convection, $\Ray$ can exceed $10^{20}$ \cite{LS23}, while experiments and simulations have been limited to at most $\Ray\approx10^{17}$ \cite{CS12}. So to estimate the heat transport in high Rayleigh regimes, we must extrapolate from scaling laws. Experiments and numerical simulations suggest that $\Nuss$ grows like a power of $\Ray$, with exponent between $0.28$ and $0.38$ \cite{LS24}, while classical heuristic arguments predict either a $\Ray^{1/3}$ scaling \cite{Mal54} or a $\Ray^{1/2}$ scaling \cite{Kra62,Spi71}. Which of the two is correct for convection between no-slip plates at finite Prandtl number remains disputed. See \cite{HFNBA12,Doe20,ISSS20,LS24}. 

Over the past several decades, a substantial body of mathematical work has established rigorous upper bounds on heat transport \cite{How63,CD96,DOR06,OS11,FGHW16,CNO16,CM25}. At finite Prandtl number, which measures the ratio of momentum diffusivity to thermal diffusivity, the unconditional estimate
\begin{equation}
    \Nuss \leq C\Ray^{1/2}
\end{equation}
holds for large $\Ray$ \cite{CD96}. In the infinite-Prandtl-number regime, where the velocity is determined instantaneously by the temperature through a Stokes equation, the corresponding upper bound improves to the $\Ray^{1/3}$ scaling. This was initially established up to logarithmic corrections \cite{CD99,DOR06,OS11}, which were subsequently removed in \cite{CM25}. These results concern no-slip boundary conditions, under which the velocity vanishes at the upper and lower boundaries.

The choice of boundary conditions plays an important role in rigorous estimates of heat transport. For the free-slip boundary conditions considered in this paper, the normal velocity and tangential stress vanish at the boundaries. In this setting, the upper bound
\begin{equation}
    \Nuss \leq C\Ray^{5/12}
\end{equation}
has been established in two dimensions for arbitrary Prandtl number \cite{WD11}, and in three dimensions in the infinite-Prandtl-number regime \cite{WD12}, with further developments for large but finite Prandtl number \cite{WW13}. Related estimates have also been obtained for Navier-slip boundary conditions, which interpolate between the no-slip and free-slip cases \cite{DNN22}.

In contrast to the extensive literature on upper bounds, much less is known about lower bounds on $\Nuss$. To our knowledge, no unconditional lower bound that grows with $\Ray$ has previously been established for the standard Rayleigh--B\'enard equations. Existing lower estimates are conditional on assumptions about the flow, such as marginal linear stability of the mean temperature profile, which have not been rigorously verified \cite{WDCK22}. Indeed, for the deterministic Rayleigh--B\'enard equations, no lower bound strictly greater than one can hold uniformly over all admissible solutions: for every Rayleigh number, the motionless conductive solution remains an admissible solution with $\Nuss=1$ \cite[Section~1]{WDCK22}. Although this solution becomes linearly unstable above a critical Rayleigh number, as first established in the classical stability analysis of Rayleigh \cite{Ray16} (see also \cite[Chapter~II]{Cha61}), linear instability alone does not provide a quantitative estimate on long-time heat transport. In particular, it does not quantify the heat transport generated by general solutions or rule out trajectories that remain close to the conductive state for arbitrarily long times. The nonlinear and infinite-dimensional nature of the dynamics makes it difficult to obtain quantitative control over such behavior. As observed in \cite[Section~1]{WDCK22}, ``[t]he concept of stability is clearly the missing ingredient, yet it is unclear how to rigorously formulate this notion as a constraint.''

A natural way to overcome this obstruction is to introduce a stochastic perturbation that acts directly on the temperature equation. Chernyshenko \cite{Che22} observed that an ever-present small noise makes an unstable solution irrelevant for long-time averages, and that, to benefit from this observation, the auxiliary functional method has to be extended to systems with noise, which has been done for finite-dimensional systems \cite{CGHP14,FGHC16}. Besides its mathematical role, such forcing has a physical motivation: random thermal fluctuations have been observed in experiments on convection, including when the fluid is at or near the conductive state \cite{WAC95,OOSA04}. In our setting, the noise is introduced exclusively in the temperature equation, with no direct stochastic forcing of the velocity. Its spatial structure is chosen so that temperature fluctuations generate non-hydrostatic buoyancy variations, which in turn induce fluid motion through the Stokes equation. Thus, the conductive equilibrium of the deterministic system is no longer an invariant configuration of the stochastically forced dynamics. Importantly, the amplitude of the stochastic forcing can be taken arbitrarily small, provided it is nonzero. The resulting lower bound does not require a positive lower threshold on the noise amplitude: however small the forcing, the estimate becomes nontrivial once the Rayleigh number is sufficiently large relative to its amplitude. In this way, persistent stochastic thermal fluctuations provide a mechanism for obtaining quantitative lower bounds on long-time heat transport.

Our main result establishes a quantitative lower bound on the long-time averaged heat transport for the stochastically forced infinite-Prandtl-number system with free-slip boundary conditions.

\begin{theorem}[Main result (informal)]
Let $0<\sigma\leq 1$ denote the amplitude of the stochastic forcing. There exist constants $c,C>0$, depending only on the domain, such that if
\begin{equation}
    \sigma^3\Ray
    \geq C\bigl(1+\log(1/\sigma)\bigr)^{192/7},
\end{equation}
then, for every deterministic initial temperature as in Theorem~\ref{thm:lower-bound}, the expected long-time averaged Nusselt number satisfies
\begin{equation}
    \underline{\Nuss}-1
    \geq c\,(\sigma^3\Ray)^{7/192}.
\end{equation}
\end{theorem}

Several features of this result are worth emphasizing. First, unlike conditional lower bounds based on assumptions about the dynamics, our estimate holds for every deterministic initial temperature in the natural $L^2$ energy class. Second, the estimate quantifies the dependence on the noise amplitude $\sigma$. In particular, for every fixed $0<\sigma\leq1$, the result implies at least polynomial growth of the long-time averaged heat transport as $\Ray\to\infty$. Thus, arbitrarily weak stochastic thermal fluctuations can produce a quantitatively nontrivial enhancement of heat transport, provided the Rayleigh number is sufficiently large. The exponent $7/192$ is not optimized: it reflects the interpolation inequalities and the choice of parameters in our argument.

Under suitable non-degeneracy assumptions on the stochastic forcing, the Markov process associated with the temperature equation admits a unique invariant probability measure and is ergodic. Such properties have been established for randomly forced infinite-Prandtl-number Rayleigh--B\'enard convection with no-slip boundary conditions \cite{FGHW16,FGR19}, see also \cite{FGHW24, FS25}. We expect that the corresponding arguments can be adapted to the free-slip setting considered here. In that case, the expected lower bound also yields an almost-sure lower bound on the long-time averaged Nusselt number for almost every stationary initial condition, by Birkhoff's ergodic theorem. We formulate this consequence in Corollary~\ref{cor:stationary}.

We focus on the infinite-Prandtl-number regime, in which viscous effects dominate inertial effects and the velocity is determined instantaneously by the temperature through a stationary Stokes equation. This approximation is relevant to highly viscous fluids, notably in models of mantle convection, and has also proved particularly useful in the mathematical analysis of heat transport \cite{Cha61,Wan04}. In this regime, the momentum equation becomes
\begin{equation}
    \begin{aligned}
        -\Delta\vu+\nabla p
            &=\Ray T\ez,\\
        \divg\vu&=0,
    \end{aligned}
\end{equation}
where $\ez$ denotes the vertical unit vector. Consequently, the velocity can be expressed as a linear function of the temperature through the Stokes solution operator. Substituting this relation into the temperature equation yields a closed evolution equation for the temperature, retaining the nonlinear advective mechanism responsible for convective heat transport while eliminating the inertial dynamics of the velocity.

The free-slip boundary conditions considered here impose
\begin{equation}
    u_2=0,\qquad
    \partial_zu_1=0
    \qquad \text{at }z=0,1,
\end{equation}
where $u_1$ and $u_2$ are the horizontal and vertical components of $\vu$, and $z$ is the vertical coordinate. These conditions are particularly convenient in the infinite-Prandtl-number formulation, since they allow the Stokes operator to be treated using its spectral structure. They also distinguish our analysis from much of the existing literature on upper bounds, which has predominantly focused on no-slip boundaries.

Our approach builds on the auxiliary functional method developed by Chernyshenko et al.~\cite{CGHP14}. The underlying idea is to establish bounds on long-time averages by constructing a suitable functional on the state space. More precisely, suppose that $F$ is an observable of a dynamical system and that $V$ is a functional bounded above along its trajectories. If
\begin{equation}
    F+\mathcal{L}V\geq B
\end{equation}
holds throughout the state space, where $\mathcal{L}$ denotes the generator of the dynamics, then, under suitable regularity and integrability assumptions, averaging along trajectories yields a lower bound of $B$ on the long-time average of $F$. For stochastic dynamics, $\mathcal{L}$ is the Markov generator and accounts for both the deterministic evolution and the contribution of the noise through It\^o's formula.

This framework provides a common perspective on several classical bounding techniques. In particular, the background method \cite{DC96,CD96} can be interpreted as a special case in which the auxiliary functional is quadratic \cite{Che22}; see also \cite{FAW22} for a review. Auxiliary functionals have been used to obtain sharp bounds for ordinary differential equations with compact invariant dynamics \cite{TGD18}, and these ideas have subsequently been extended to the two-dimensional Navier--Stokes equations \cite{RT22}. For polynomial systems, the search for suitable functionals can also be formulated as a sum-of-squares optimization problem, including in the presence of stochastic forcing \cite{CGHP14,FGHC16}.

In this work, we extend the auxiliary functional framework to Markov processes and use it to obtain lower bounds on long-time averages in stochastically forced infinite-Prandtl-number convection. A key feature of our approach is a relaxation of the usual pointwise auxiliary functional inequality. Rather than requiring the desired inequality to hold throughout the state space, we divide the state space into two regions. On the first region, the inequality holds directly. On the complementary region, where the instantaneous contribution is insufficient, we exploit the subsequent evolution of the dynamics over a controlled finite time interval. By averaging the contribution along these short trajectories, we recover a positive contribution that compensates for the failure of the pointwise inequality. This yields a finite-time averaged version of the auxiliary functional principle, providing additional flexibility in the construction of functionals for nonlinear stochastic evolution equations. We develop this framework in Section~\ref{sec:window} and apply it to the Nusselt number to establish the lower bound in Theorem~\ref{thm:lower-bound}. Related conditions on averages over a time interval appear in the stability theory of Markov processes \cite{MD94,MT94}.

In this application, at the states where the pointwise inequality fails, we project the solution onto an unstable eigenfunction of the linearization around the initial horizontal mean temperature profile, frozen in time. We take this projection in the inner product defined by the convective heat flux: the linearization is symmetric in this inner product, and a multiple of the square of the projection bounds the heat flux from below.

\subsection{The main result}\label{subsec:statement}

We consider the infinite-Prandtl-number Rayleigh--B\'enard Boussinesq system on $D\coloneqq[0,L]\times[0,1]$:
\begin{equation}\label{eq:T-system}
\begin{aligned}
    &{-}\Delta\vu+\nabla p=\Ray\,T\ez,\\
    &\divg\vu=0,\\
    &dT+\vu\cdot\nabla T\,dt
        =\Delta T\,dt+\sigma\sum_{n=1}^\infty\zeta_n\,dW^n_t.
\end{aligned}
\end{equation}
Here, $\vu:D\times[0,\infty)\to\R^2$ is the velocity field, $T:D\times[0,\infty)\to\R$ is the temperature, and $p:D\times[0,\infty)\to\R$ is the pressure. The system is parametrized by the Rayleigh number $\Ray>0$.

We consider solutions that are periodic in $x$ and satisfy the free-slip boundary conditions
\begin{equation}\label{eq:T-bc}
\begin{aligned}
    &T=1\quad\text{at }z=0,\qquad T=0\quad\text{at }z=1,\\
    &u_2=\partial_z u_1=0\quad\text{at }z=0,1.
\end{aligned}
\end{equation}
We also impose the zero-mean condition
\begin{equation}\label{eq:zero-mean}
    \int_0^1\!\int_0^L u_1\,dx\,dz=0,
\end{equation}
which fixes the horizontally constant component of the velocity, left undetermined by the free-slip boundary conditions.

The stochastic forcing in \eqref{eq:T-system} is given by
\begin{equation}\label{eq:noise}
    \zeta_n(x,z)\coloneqq
    C_L\,\frac{k_0}{k_0^2+n^2\pi^2}
    \cos(k_0x)\sin(n\pi z),
\end{equation}
where $k_0\coloneqq2\pi/L$, the $W^n$ are independent standard Brownian motions, and $0<\sigma\leq1$ is the noise amplitude. The constant $C_L>0$, depending only on $L$, is chosen so that
\begin{equation}\label{eq:noise-variance}
    \sum_{n=1}^\infty\norm{\zeta_n}_{L^2(D)}^2=1.
\end{equation}
Thus, the covariance operator of the $L^2(D)$-valued noise has trace $\sigma^2$ per unit time.

The instantaneous spatially averaged vertical heat transport is
\begin{equation}\label{eq:instantaneous-nusselt}
    \mathcal N(T)
    \coloneqq
    \frac1L\int_0^1\!\int_0^L
    \bigl(u_2T-\partial_zT\bigr)\,dx\,dz.
\end{equation}
Here, $\vu$ is determined by $T$ through the Stokes equation in \eqref{eq:T-system}. This expression measures the vertical heat transport relative to the transport by pure conduction, for which $\mathcal N(1-z)=1$.

We define the expected long-time averaged Nusselt number by
\begin{equation}\label{eq:nusselt}
    \underline{\Nuss}
    \coloneqq
    \liminf_{t\to\infty}
    \frac1t\,\E_{\theta_0}\int_0^t\mathcal N(T(s))\,ds,
\end{equation}
where $\E_{\theta_0}$ denotes expectation with respect to the stochastic dynamics initialized at $T(0)=1-z+\theta_0$, for a deterministic perturbation $\theta_0\in L^2(D)$. We use a liminf since the existence of the long-time limit is not established, and our objective is to obtain a lower bound without requiring its existence.

\begin{theorem}\label{thm:lower-bound}
Let the noise be as in \eqref{eq:noise}, with $0<\sigma\leq1$. There exist constants $c,C>0$, depending only on $L$, such that if
\begin{equation}\label{eq:parameter-regime}
    \sigma^3\Ray
    \geq C\bigl(1+\log(1/\sigma)\bigr)^{192/7},
\end{equation}
then, for every deterministic initial temperature $T(0)$ satisfying
\[
    T(0)-(1-z)\in L^2(D),
\]
we have
\begin{equation}\label{eq:main-bound}
    \underline{\Nuss}-1
    \geq c\,(\sigma^3\Ray)^{7/192}.
\end{equation}
\end{theorem}

The proof of Theorem~\ref{thm:lower-bound} is given in Section~\ref{sec:proof-main}.

Under additional assumptions on the stochastic forcing, the Markov process associated with the temperature equation admits a unique invariant probability measure and is ergodic. In this case, the expected lower bound in Theorem~\ref{thm:lower-bound} also gives an almost-sure lower bound for the long-time averaged heat transport of stationary solutions.

\begin{corollary}\label{cor:stationary}
Let $\sigma$ and $\Ray$ satisfy the assumptions of Theorem~\ref{thm:lower-bound}, with the same constants $c$ and $C$. Assume, in addition, that the associated Markov semigroup admits a unique invariant probability measure $\mu$ on $L^2(D)$ and that the corresponding stationary Markov process is ergodic. Assume also that $\theta\mapsto\mathcal N(1-z+\theta)$ is integrable with respect to $\mu$. Then the following hold:
\begin{enumerate}[label=(\roman*)]
    \item The stationary Nusselt number
    \begin{equation*}
        \Nuss_\mu
        \coloneqq
        \int_{L^2(D)}
        \mathcal N(1-z+\theta)\,\mu(d\theta)
    \end{equation*}
    satisfies
    \begin{equation*}
        \Nuss_\mu-1
        \geq c\,(\sigma^3\Ray)^{7/192}.
    \end{equation*}

    \item For $\mu$-almost every initial perturbation
    $\theta_0\in L^2(D)$, almost surely,
    \begin{equation*}
        \lim_{t\to\infty}
        \frac1t\int_0^t\mathcal N(T(s))\,ds-1
        \geq c\,(\sigma^3\Ray)^{7/192},
    \end{equation*}
    where $T(0)=1-z+\theta_0$.
\end{enumerate}
\end{corollary}

Part (i) follows by integrating~\eqref{eq:main-bound} with respect to $\mu$, using the invariance of $\mu$ and Fatou's lemma. Part (ii) follows from part (i) and Birkhoff's ergodic theorem.

\begin{remark}
The unique ergodicity assumption can be established, for suitable forcing, using the general strategy developed in \cite{FGHW16,FGR19}, where uniqueness of the invariant measure is proved for randomly forced Rayleigh--B\'enard convection in the infinite-Prandtl-number regime. Although these works consider no-slip boundary conditions, we expect their arguments to extend to the free-slip setting considered here. In particular, the free-slip Stokes operator admits a convenient spectral representation that facilitates the corresponding analysis. The forcing in these works also acts on the horizontal mean of the temperature. We expect that the proof of Theorem~\ref{thm:lower-bound} can be adapted to such a forcing (compare Remark~\ref{rem:general-noise}). We do not pursue this question here, and instead assume that the stochastic forcing satisfies suitable non-degeneracy conditions guaranteeing unique ergodicity.
\end{remark}

\begin{remark}\label{rem:general-noise}
The proof of Theorem~\ref{thm:lower-bound} can be carried out without any loss if the noise in \eqref{eq:T-system} is
\[
    \sigma\sum_{n=1}^\infty\zeta_n\,dW^n_t+\sum_{m=1}^\infty\chi_m\,d\widetilde W^m_t,
\]
where the $\widetilde W^m$ are standard Brownian motions, all the $W^n$ and $\widetilde W^m$ are independent, and the functions $\chi_m\in L^2(D)$ satisfy
\[
    \int_0^L\chi_m(x,z)\,dx=0\quad\text{for almost every }z\in(0,1),\qquad\sigma^2+\sum_{m=1}^\infty\norm{\chi_m}_{L^2(D)}^2\leq1.
\]
Indeed, the constant $\sigma^2$ in \eqref{eq:energy-ito} becomes $\sigma^2+\sum_m\norm{\chi_m}_{L^2(D)}^2\leq1$, and, since the Brownian motions are independent, $q$ in \eqref{eq:q-formula} gains the nonnegative term $\sum_m(\chi_m,K\Psi)^2$. The proof uses that the noise has zero horizontal mean. It does not cover noise acting on finitely many modes, or noise in the horizontal mean.
\end{remark}

The remainder of the paper is organized as follows. Section~\ref{sec:window} develops the finite-time auxiliary functional framework for Markov processes and proves Lemma~\ref{lem:window}, which is central to the proof of Theorem~\ref{thm:lower-bound}. Section~\ref{sec:system} formulates the system in terms of the perturbation $\theta=T-(1-z)$ and establishes the estimates needed for the main result. Section~\ref{sec:proof-main} contains the proof of Theorem~\ref{thm:lower-bound}. The construction of the one-dimensional projection used in the proof and the choice of parameters are carried out in Sections~\ref{sec:coordinate} and~\ref{sec:verify}, respectively. Section~\ref{sec:proofs} contains the remaining technical arguments, and Appendix~\ref{sec:wellposed} establishes well-posedness and the Markov property.

\subsection{Conventions and notation}\label{subsec:conventions}

Functions on $D$ are assumed to be $L$-periodic in $x$, so that we identify them with functions on the cylinder $(\R/L\mathbb Z)\times[0,1]$. We use the inner product and norm
\[
    (f,g)\coloneqq\int_0^1\!\int_0^L fg\,dx\,dz,
    \qquad
    \norm{f}_{L^2(D)}\coloneqq(f,f)^{1/2},
    \qquad f,g\in L^2(D).
\]
For vector fields $\vu,\vv:D\to\R^2$, we write
\[
    (\vu,\vv)\coloneqq(u_1,v_1)+(u_2,v_2)
\]
and
\[
    (\nabla\vu,\nabla\vv)
    \coloneqq\sum_{i,j}(\partial_j u_i,\partial_j v_i).
\]

The following functions form an orthonormal basis of $L^2(D)$:
\begin{equation}\label{eq:fourier-basis}
\begin{aligned}
    e_n&\coloneqq\frac{\sqrt2}{L^{1/2}}\sin(n\pi z),\\
    c_{k,n}&\coloneqq\frac{2}{L^{1/2}}\cos(kx)\sin(n\pi z),\\
    s_{k,n}&\coloneqq\frac{2}{L^{1/2}}\sin(kx)\sin(n\pi z),
\end{aligned}
\end{equation}
where $n\geq1$ and $k\in\{k_0,2k_0,\dots\}$, with $k_0=2\pi/L$. The functions $e_n$ depend only on $z$ and correspond to the horizontally constant modes. Thus, every $f\in L^2(D)$ admits the expansion
\[
    f
    =
    \sum_{n\geq1}(f,e_n)e_n
    +
    \sum_{\substack{k\in\{k_0,2k_0,\dots\}\\n\geq1}}
    \Bigl(
        (f,c_{k,n})c_{k,n}
        +(f,s_{k,n})s_{k,n}
    \Bigr).
\]

We denote by $H^1_0(D)$ the space of functions in $H^1(D)$ that are $L$-periodic in $x$ and have zero trace at $z=0,1$. In particular, the notation $H^1_0(D)$ here refers only to the vertical boundary conditions, with periodicity imposed in the horizontal direction. We equip this space with the gradient norm
\begin{equation}\label{eq:H1-fourier}
\begin{aligned}
    \norm{\theta}_{H^1_0(D)}^2
    \coloneqq{}&
    \sum_{\substack{k\in\{k_0,2k_0,\dots\}\\n\geq1}}
    (k^2+n^2\pi^2)
    \bigl[
        (\theta,c_{k,n})^2
        +(\theta,s_{k,n})^2
    \bigr]\\
    &+\sum_{n\geq1}n^2\pi^2(\theta,e_n)^2.
\end{aligned}
\end{equation}
This norm is equivalent to the standard $H^1$ norm on $H^1_0(D)$ by Poincar\'e's inequality. For $\theta\in H^1_0(D)$, integration by parts gives
\[
    \norm{\nabla\theta}_{L^2(D)}
    =\norm{\theta}_{H^1_0(D)}.
\]

We write $H^{-1}(D)$ for the dual of $H^1_0(D)$, with the dual norm induced by the gradient norm. Each $f\in L^2(D)$ defines an element of $H^{-1}(D)$ through the duality pairing
\[
    \langle f,\vartheta\rangle_{H^{-1}\times H^1_0}=(f,\vartheta),
    \qquad \vartheta\in H^1_0(D).
\]
The corresponding norm is given by
\begin{equation}\label{eq:Hminus1-fourier}
\begin{aligned}
    \norm{f}_{H^{-1}(D)}^2
    =
    \sum_{\substack{k\in\{k_0,2k_0,\dots\}\\n\geq1}}
    \frac{(f,c_{k,n})^2+(f,s_{k,n})^2}
    {k^2+n^2\pi^2}
    +\sum_{n\geq1}\frac{(f,e_n)^2}{n^2\pi^2}.
\end{aligned}
\end{equation}

We write $\mathbb P_x$ and $\E_x$ for the probability and expectation associated with the process started at $x$. In Section~\ref{sec:window}, the process is a general Markov process $X_t$ with values in a state space $\mathcal S$, started at $X_0=x\in\mathcal S$. In Sections~\ref{sec:system}--\ref{sec:proofs} and Appendix~\ref{sec:wellposed}, the process is the perturbation $\theta=T-(1-z)$, with $\theta(0)=\theta_0\in L^2(D)$, solving \eqref{eq:theta-equation} in Section~\ref{subsec:solutions}.

Throughout the estimates in Sections~\ref{sec:system}--\ref{sec:proofs} and Appendix~\ref{sec:wellposed}, $C$ denotes a sufficiently large positive constant and $c$ a sufficiently small positive constant. Their values may change from line to line, but they depend only on $L$.

\section{The auxiliary functional method for Markov processes}\label{sec:window}

Throughout this section, we fix a time-homogeneous Markov process $X_t$ on a filtered probability space $(\Omega, \mathcal{F}, \mathcal{F}_{t \geq 0}, \mathbb{P})$, with values in a state space $\mathcal S$ and with transition kernel $P_t(x,A)\coloneqq\mathbb P_x(X_t\in A)$, for $t\geq0$, $x\in\mathcal S$ and Borel sets $A\subset\mathcal S$. Let $F:\mathcal S\to\R$ be the quantity whose long-time average we want to estimate from below, that is, we are interested in
\[
   \liminf_{S\to\infty}\frac1S\,\E_x\int_0^SF(X_t)\,dt.
\]
We also assume that $(t,\omega)\mapsto X_t(\omega)$, $[0,\infty)\times\Omega\to\mathcal S$, is jointly measurable and that $F$ is a Borel function, to allow time integrals such as $\int_0^SF(X_t)\,dt$ to be random variables, and to allow exchanging $\E_x$ and $\int_0^S$.

\subsection{Extension of the method}\label{subsec:extension}

We present an extension of the auxiliary functional argument of Chernyshenko et al.\ \cite{CGHP14,Che22} to a stochastic setting (see \cite{CGHP14,FGHC16} for stationary expectations of finite-dimensional stochastic systems), which allows functionals $V$ that may be unbounded above and not regular enough for the generator to exist. We let $\LL V$ refer to an extended generator of $V$:

\begin{definition}[{\cite[Section~1.3]{MT93}}]\label{def:extended-generator}
A Borel function $\LL V:\mathcal S\to\R$ is an \emph{extended generator} of $V:\mathcal S\to\R$ with respect to the process $X_t$ if, for every $x\in\mathcal S$ and $S>0$,
\begin{equation}\label{eq:LV-integrable}
    \E_x\int_0^S\abs{\LL V(X_t)}\,dt<\infty,
\end{equation}
and Dynkin's formula holds:
\begin{equation}\label{eq:dynkin}
    \E_xV(X_S)-V(x)=\E_x\int_0^S\LL V(X_t)\,dt.
\end{equation}
\end{definition}

For example, an extended generator is usually found as the drift in It\^o's formula. If, for every $x\in\mathcal S$, $\mathbb P_x$-almost surely,
\begin{equation}\label{eq:ito-drift}
    V(X_t)=V(x)+\int_0^t\LL V(X_s)\,ds+M_t\qquad\text{for }t\geq0,
\end{equation}
where $\LL V:\mathcal S\to\R$ is Borel and satisfies \eqref{eq:LV-integrable} for every $S>0$, and $M$ is a martingale with $M_0=0$, then taking expectations in \eqref{eq:ito-drift} yields Dynkin's formula, and therefore $\LL V$ is an extended generator of $V$.

Let $V:\mathcal S\to\R$ have an extended generator $\LL V$, and let $B\in\R$ and $x\in\mathcal S$ with
\begin{equation}\label{eq:F-integrable}
    \E_x\int_0^S\abs{F(X_t)}\,dt<\infty\qquad\text{for every }S>0.
\end{equation}
If
\begin{equation}\label{eq:af-conditions}
\begin{aligned}
    &\text{(i)}\quad\liminf_{S\to\infty}\E_x\frac1S\int_0^S(F+\LL V)(X_t)\,dt\geq B,\\
    &\text{(ii)}\quad\liminf_{S\to\infty}\E_x\frac1S\int_0^SF(X_t)\,dt\geq\liminf_{S\to\infty}\E_x\frac1S\int_0^S(F+\LL V)(X_t)\,dt,
\end{aligned}
\end{equation}
then,
\begin{equation}\label{eq:af-bound}
    \liminf_{S\to\infty}\E_x\frac1S\int_0^SF(X_t)\,dt\geq\liminf_{S\to\infty}\E_x\frac1S\int_0^S(F+\LL V)(X_t)\,dt\geq B.
\end{equation}

\begin{remark}\label{rem:upper-bounds}
The implication from \eqref{eq:af-conditions} to \eqref{eq:af-bound} also gives upper bounds: replacing $F$, $V$ and $B$ by $-F$, $-V$ and $-B$ shows that, if conditions (i) and (ii) hold with $\leq$ in place of $\geq$ and with $\limsup$ in place of $\liminf$, then
\[
    \limsup_{S\to\infty}\E_x\frac1S\int_0^SF(X_t)\,dt\leq B.
\]
\end{remark}

\begin{remark}\label{rem:af-method}
The auxiliary functional method in the form \eqref{eq:af-conditions}, \eqref{eq:af-bound} recovers several known bounds. By \eqref{eq:dynkin},
\[
    \E_x\frac1S\int_0^S(F+\LL V)(X_t)\,dt=\E_x\frac1S\int_0^SF(X_t)\,dt+\frac{\E_xV(X_S)-V(x)}S.
\]
Therefore, if $V$ is bounded above and $F+\LL V\geq B$ on $\mathcal S$, conditions (i) and (ii) hold. In this case, the bound \eqref{eq:af-bound} is the auxiliary functional method of Chernyshenko et al.\ \cite{CGHP14,Che22} for deterministic systems, where the derivative of $V$ along a solution takes the place of $\LL V$. For Markov processes, this bound is a form of the comparison theorem of Meyn and Tweedie \cite[Theorem~1.1]{MT93}. For stationary distributions of Markov processes, the analogous upper bound is given by Glynn and Zeevi \cite[Theorem~1]{GZ08}.
\end{remark}

\subsection{The short-time average condition}\label{subsec:short-time}

We start with an auxiliary functional $V_0:\mathcal S\to\R$, bounded above, with an extended generator $\LL V_0$, such that, for a Borel function $Q:\mathcal S\to[0,\infty)$ and a constant $c_0\geq0$,
\begin{equation}\label{eq:balance}
    F+\LL V_0\geq Q-c_0 .
\end{equation}

By itself $V_0$ only ensures $F+\LL V_0\geq-c_0$, and hence conditions (i) and (ii) in \eqref{eq:af-conditions} with $V=V_0$ yield only the lower bound $-c_0$ (see Remark~\ref{rem:af-method}).

To improve the lower bound $-c_0$, fix $\beta>0$. Where $Q\geq\beta$ we already have $F+\LL V_0\geq Q-c_0\geq\beta-c_0$, and therefore it remains to treat the set
\begin{equation}\label{eq:low-set-def}
    \mathcal S_{Q<\beta}\coloneqq\{x\in\mathcal S:Q(x)<\beta\}.
\end{equation}

On $\mathcal S_{Q<\beta}$ we only know $Q(x)<\beta$, and consequently we replace $Q(x)$ by the expected values $\E_xQ(X_t)$ at later times $t>0$: we require that, for some $\tau>0$,
\begin{equation}\label{eq:window-condition}
    \E_x\,\frac1\tau\int_0^\tau Q(X_t)\,dt\geq\beta\qquad\text{for every }x\in\mathcal S_{Q<\beta}.
\end{equation}
The following lemma shows that the condition \eqref{eq:window-condition} improves the lower bound $-c_0$ to $\frac\beta2-c_0$.

\begin{lemma}\label{lem:window}
Let $X_t$ have a transition kernel that is jointly measurable in $(t,x)$, and let $F\geq0$ satisfy \eqref{eq:F-integrable} for every $x\in\mathcal S$, and let $V_0:\mathcal S\to\R$ be bounded above and have an extended generator $\LL V_0$. Let $Q:\mathcal S\to[0,\infty)$ be Borel and $c_0\geq0$ be such that \eqref{eq:balance} holds, that is, $F+\LL V_0\geq Q-c_0$ on $\mathcal S$. Fix $\beta>0$, and suppose \eqref{eq:window-condition} holds for some $\tau>0$ independent of $x\in\mathcal S_{Q<\beta}=\{x\in\mathcal S:Q(x)<\beta\}$. Then for every $x\in\mathcal S$,
\begin{equation}\label{eq:window-bound}
    \liminf_{S\to\infty}\frac1S\,\E_x\int_0^S F(X_t)\,dt\geq\frac\beta2-c_0 .
\end{equation}
\end{lemma}

\begin{proof}
The proof is inspired by Meyn and Down \cite[Lemma~2.3]{MD94}. Fix $x\in\mathcal S$. Since $Q\geq0$, the condition \eqref{eq:window-condition} implies
\begin{equation}\label{eq:window-everywhere}
    Q(y)+\E_y\frac1\tau\int_0^\tau Q(X_r)\,dr\geq\beta\qquad\text{for every }y\in\mathcal S,
\end{equation}
because the first term on the left side is at least $\beta$ if $Q(y)\geq\beta$, and the second term is at least $\beta$ if $Q(y)<\beta$.

We evaluate \eqref{eq:window-everywhere} at the states $y=X_t$. By the Markov property, $\E_x\E_{X_t}Q(X_r)=\E_xQ(X_{t+r})$. Since the transition kernel $P_t(y,A)$ is jointly measurable in $(t,y)$, the map
\[
    (t,y)\mapsto\E_yQ(X_t)=\int_{\mathcal S}Q(z)\,P_t(y,dz)
\]
is also jointly measurable. Indeed, this holds if $Q$ is the indicator function of a Borel set, and hence for every Borel $Q\geq0$ by the monotone convergence theorem. Therefore, Tonelli's theorem applies, and taking $\E_x$ in \eqref{eq:window-everywhere} with $y=X_t$, we obtain
\begin{equation}\label{eq:window-at-t}
    \E_xQ(X_t)+\frac1\tau\int_0^\tau\E_xQ(X_{t+r})\,dr\geq\beta\qquad\text{for every }t\geq0.
\end{equation}

For $S>0$, averaging \eqref{eq:window-at-t} over $t\in[0,S]$ leads to
\begin{equation}\label{eq:Q-average-bound}
\begin{aligned}
    \beta&\leq\frac1S\int_0^S\E_xQ(X_t)\,dt+\frac1{S\tau}\int_0^\tau\Bigl\{\int_0^S\E_xQ(X_{t+r})\,dt\Bigr\}\,dr\\
    &=\frac1S\int_0^S\E_xQ(X_t)\,dt+\frac1{S\tau}\int_0^\tau\Bigl\{\int_r^{S+r}\E_xQ(X_t)\,dt\Bigr\}\,dr\\
    &\leq\frac2S\int_0^{S+\tau}\E_xQ(X_t)\,dt,
\end{aligned}
\end{equation}
where the last inequality holds because $Q\geq0$ and $r\leq\tau$.

Finally, we pass from $Q$ to $F$. By \eqref{eq:balance}, then by Dynkin's formula \eqref{eq:dynkin} for $V_0$, and then by $V_0\leq\sup_{\mathcal S}V_0$,
\[
\begin{aligned}
    \int_0^{S+\tau}\E_xQ(X_t)\,dt&\leq\E_x\int_0^{S+\tau}\bigl(F+\LL V_0+c_0\bigr)(X_t)\,dt\\
    &=\E_x\int_0^{S+\tau}F(X_t)\,dt+\E_xV_0(X_{S+\tau})-V_0(x)+c_0(S+\tau)\\
    &\leq\E_x\int_0^{S+\tau}F(X_t)\,dt+\sup_{\mathcal S}V_0-V_0(x)+c_0(S+\tau).
\end{aligned}
\]
Then \eqref{eq:Q-average-bound} implies
\[
    \frac\beta2\leq\frac{S+\tau}S\Bigl(\E_x\frac1{S+\tau}\int_0^{S+\tau}F(X_t)\,dt+c_0\Bigr)+\frac{\sup_{\mathcal S}V_0-V_0(x)}S,
\]
and letting $S\to\infty$ yields \eqref{eq:window-bound}.
\end{proof}

\begin{remark}\label{rem:spirit}
Condition \eqref{eq:window-condition} is in the spirit of Meyn and Down \cite[Lemma~2.5]{MD94} and, in discrete time, Meyn and Tweedie \cite[Section~4(c)]{MT94}, which estimate the average of $\E_xV(X_t)-V(x)$ over a time interval in place of $\LL V$.
\end{remark}

\subsection{The short-time average functional}\label{subsec:window-functional}

The lower bound \eqref{eq:window-bound} of Lemma~\ref{lem:window} can also be obtained from conditions (i) and (ii) in \eqref{eq:af-conditions}, with an explicit functional $V$. For $\tau>0$, the \emph{short-time average functional} $G_\tau:\mathcal S\to[0,\infty]$ is
\begin{equation}\label{eq:window-functional}
    G_\tau(x)\coloneqq\E_x\int_0^\tau\Bigl(1-\frac t\tau\Bigr)Q(X_t)\,dt .
\end{equation}

\begin{lemma}\label{lem:window-generator}
Let the transition kernel $P_t(x,A)$ of $X_t$ be jointly measurable in $(t,x)$ for every Borel set $A\subset\mathcal S$. Let $Q:\mathcal S\to[0,\infty)$ be Borel, with
\begin{equation}\label{eq:Q-integrable}
    \E_x\int_0^SQ(X_t)\,dt<\infty\qquad\text{for every $x\in\mathcal S$ and $S>0$.}
\end{equation}
Then, for every $\tau>0$, $G_\tau$ in \eqref{eq:window-functional} is finite and Borel, and $G_\tau$ has an extended generator, given by
\begin{equation}\label{eq:window-generator}
    \LL G_\tau(x)=\frac1\tau\int_0^\tau\bigl(\E_xQ(X_r)-Q(x)\bigr)\,dr=-Q(x)+\E_x\frac1\tau\int_0^\tau Q(X_t)\,dt.
\end{equation}
\end{lemma}

\begin{proof}
Define $\LL G_\tau$ by \eqref{eq:window-generator}. Since the transition kernel $P_t(x,A)$ is jointly measurable in $(t,x)$, so is $(t,x)\mapsto\E_xQ(X_t)$, and therefore $G_\tau$ and $\LL G_\tau$ are Borel, by Tonelli's theorem. Moreover, by \eqref{eq:Q-integrable},
\[
    G_\tau(x)\leq\E_x\int_0^\tau Q(X_t)\,dt<\infty.
\]
By the Markov property, $\E_x\E_{X_t}Q(X_r)=\E_xQ(X_{t+r})$, and consequently, by Tonelli's theorem, for $0\leq r\leq\tau$,
\begin{equation}\label{eq:G-shift}
    \E_x\int_0^S\E_{X_t}Q(X_r)\,dt=\int_0^S\E_xQ(X_{t+r})\,dt=\E_x\int_r^{S+r}Q(X_t)\,dt\leq\E_x\int_0^{S+\tau}Q(X_t)\,dt.
\end{equation}
The definition \eqref{eq:window-generator} of $\LL G_\tau$ implies
\[
    \abs{\LL G_\tau(X_t)}\leq\frac1\tau\int_0^\tau\E_{X_t}Q(X_r)\,dr+Q(X_t),
\]
and hence, averaging \eqref{eq:G-shift} over $r\in[0,\tau]$,
\[
\begin{aligned}
    \E_x\int_0^S\abs{\LL G_\tau(X_t)}\,dt&\leq\frac1\tau\int_0^\tau\E_x\int_0^S\E_{X_t}Q(X_r)\,dt\,dr+\E_x\int_0^SQ(X_t)\,dt\\
    &\leq\E_x\int_0^{S+\tau}Q(X_t)\,dt+\E_x\int_0^SQ(X_t)\,dt<\infty.
\end{aligned}
\]

Finally, by \eqref{eq:window-functional}, the Markov property and Fubini's theorem, the left side of Dynkin's formula \eqref{eq:dynkin} for $G_\tau$ is
\begin{equation}\label{eq:G-left}
    \E_xG_\tau(X_S)-G_\tau(x)=\int_0^\tau\Bigl(1-\frac t\tau\Bigr)\bigl(\E_xQ(X_{S+t})-\E_xQ(X_t)\bigr)\,dt.
\end{equation}
By \eqref{eq:window-generator} and the equalities in \eqref{eq:G-shift}, the right side of \eqref{eq:dynkin} for $G_\tau$ is
\begin{equation}\label{eq:G-right}
\begin{aligned}
    \E_x\int_0^S\LL G_\tau(X_t)\,dt&=\frac1\tau\int_0^\tau\Bigl(\E_x\int_0^S\E_{X_t}Q(X_r)\,dt-\E_x\int_0^SQ(X_t)\,dt\Bigr)\,dr\\
    &=\frac1\tau\int_0^\tau\Bigl(\int_r^{S+r}\E_xQ(X_t)\,dt-\int_0^S\E_xQ(X_t)\,dt\Bigr)\,dr\\
    &=\frac1\tau\int_0^\tau\Bigl(\int_S^{S+r}\E_xQ(X_t)\,dt-\int_0^r\E_xQ(X_t)\,dt\Bigr)\,dr\\
    &=\frac1\tau\int_0^\tau\int_0^r\bigl(\E_xQ(X_{S+t})-\E_xQ(X_t)\bigr)\,dt\,dr\\
    &=\int_0^\tau\Bigl(1-\frac t\tau\Bigr)\bigl(\E_xQ(X_{S+t})-\E_xQ(X_t)\bigr)\,dt\\
    &=\E_xG_\tau(X_S)-G_\tau(x),
\end{aligned}
\end{equation}
where the next-to-last equality follows from Fubini's theorem, since $\frac1\tau\int_t^\tau dr=1-\frac t\tau$, and the last equality is \eqref{eq:G-left}. Hence, \eqref{eq:dynkin} holds for $G_\tau$.
\end{proof}

Let $F$, $V_0$, $Q$, $c_0$, $\beta$ and $\tau$ be as in Lemma~\ref{lem:window}. Adding \eqref{eq:balance} and one half of \eqref{eq:window-generator} yields
\begin{equation}\label{eq:window-certificate}
    F+\LL\bigl(V_0+\tfrac12G_\tau\bigr)\geq\frac12\Bigl(Q(x)+\E_x\frac1\tau\int_0^\tau Q(X_t)\,dt\Bigr)-c_0 .
\end{equation}
Adding all of $G_\tau$ would cancel $Q(x)$, and then we would need a lower bound on the left side of \eqref{eq:window-condition} at every state. The factor $\frac12$ keeps half of $Q(x)$ in \eqref{eq:window-certificate}.

\begin{proposition}\label{prop:window-functional}
Under the assumptions of Lemma~\ref{lem:window}, the functional $V\coloneqq V_0+\frac12G_\tau$ has an extended generator $\LL V\coloneqq\LL V_0+\tfrac12\LL G_\tau$ with
\begin{equation}\label{eq:window-pointwise}
    F+\LL V\geq\frac\beta2-c_0\qquad\text{on }\mathcal S,
\end{equation}
and, for every $x\in\mathcal S$, conditions (i) and (ii) in \eqref{eq:af-conditions} hold with this $V$ and with $B=\frac\beta2-c_0$. In particular, \eqref{eq:af-bound} yields the lower bound \eqref{eq:window-bound} of Lemma~\ref{lem:window}.
\end{proposition}

\begin{proof}
We show, for every $x\in\mathcal S$, conditions (i) and (ii) in \eqref{eq:af-conditions} with $V\coloneqq V_0+\frac12G_\tau$, where $G_\tau$ is as in \eqref{eq:window-functional}, and with $B=\frac\beta2-c_0$.
Since \eqref{eq:balance} implies $Q\leq F+\LL V_0+c_0$, we have
\[
    \E_x\int_0^SQ(X_t)\,dt\leq\E_x\int_0^S\bigl(F+\abs{\LL V_0}+c_0\bigr)(X_t)\,dt<\infty,
\]
where the right side is finite by \eqref{eq:F-integrable}, and by \eqref{eq:LV-integrable} for $\LL V_0$. Therefore, \eqref{eq:Q-integrable} holds, and Lemma~\ref{lem:window-generator} applies to $G_\tau$. Consequently, $\LL V=\LL V_0+\frac12\LL G_\tau$ is an extended generator of $V$.

\emph{Condition (i).} Since $Q\geq0$, the right side of \eqref{eq:window-certificate} is at least $\frac\beta2-c_0$ where $Q\geq\beta$, and also in $\mathcal S_{Q<\beta}$ by \eqref{eq:window-condition}. Therefore,
\[
    F+\LL V\geq\frac\beta2-c_0\qquad\text{on }\mathcal S,
\]
and consequently condition (i) holds.

\emph{Condition (ii).} Since $G_\tau$ is not known to be bounded above, we estimate the term $\E_xG_\tau(X_S)$ directly. Evaluating \eqref{eq:window-functional} at the state $X_S$ yields
\[
    G_\tau(X_S)=\int_0^\tau\Bigl(1-\frac t\tau\Bigr)\E_{X_S}Q(X_t)\,dt.
\]
Taking $\E_x$, and using the Markov property and Tonelli's theorem, we obtain
\[
    \E_xG_\tau(X_S)=\int_0^\tau\Bigl(1-\frac t\tau\Bigr)\E_xQ(X_{S+t})\,dt.
\]
Since $0\leq1-\frac t\tau\leq1$ on $[0,\tau]$ and $Q\geq0$, substituting $r=S+t$ yields
\begin{equation}\label{eq:G-at-XS}
    \E_xG_\tau(X_S)\leq\int_0^\tau\E_xQ(X_{S+t})\,dt=\int_S^{S+\tau}\E_xQ(X_r)\,dr.
\end{equation}

Subtracting Dynkin's formula \eqref{eq:dynkin} for $V_0$ at the time $S$ from \eqref{eq:dynkin} for $V_0$ at the time $S+\tau$ leads to
\[
    \E_x\int_S^{S+\tau}\LL V_0(X_r)\,dr=\E_xV_0(X_{S+\tau})-\E_xV_0(X_S).
\]
Thus, integrating \eqref{eq:balance} over $[S,S+\tau]$, we obtain
\[
    \int_S^{S+\tau}\E_xQ(X_r)\,dr\leq\E_x\int_S^{S+\tau}F(X_r)\,dr+c_0\tau+\E_xV_0(X_{S+\tau})-\E_xV_0(X_S).
\]
Then \eqref{eq:G-at-XS} implies
\begin{equation}\label{eq:G-control}
    \frac12\E_xG_\tau(X_S)\leq\frac12\E_x\int_S^{S+\tau}F(X_r)\,dr+\frac{c_0\tau}2+\frac12\E_xV_0(X_{S+\tau})-\frac12\E_xV_0(X_S).
\end{equation}

By \eqref{eq:dynkin} and the definition $V=V_0+\frac12G_\tau$,
\[
    \E_x\int_0^S\LL V(X_t)\,dt=\E_xV(X_S)-V(x)=\E_xV_0(X_S)+\frac12\E_xG_\tau(X_S)-V(x).
\]
Therefore, by \eqref{eq:G-control} and $V_0\leq\sup_{\mathcal S}V_0$, and then because $F\geq0$,
\[
\begin{aligned}
    \E_x\int_0^S(F+\LL V)(X_t)\,dt&=\E_x\int_0^SF(X_t)\,dt+\E_xV_0(X_S)+\frac12\E_xG_\tau(X_S)-V(x)\\
    &\leq\E_x\int_0^SF(X_t)\,dt+\frac12\E_x\int_S^{S+\tau}F(X_r)\,dr+\sup_{\mathcal S}V_0+\frac{c_0\tau}2-V(x)\\
    &\leq\E_x\int_0^{S+\tau}F(X_t)\,dt+\sup_{\mathcal S}V_0+\frac{c_0\tau}2-V(x).
\end{aligned}
\]
Dividing by $S$, we obtain
\[
    \E_x\frac1S\int_0^S(F+\LL V)(X_t)\,dt\leq\frac{S+\tau}S\,\E_x\frac1{S+\tau}\int_0^{S+\tau}F(X_t)\,dt+\frac{\sup_{\mathcal S}V_0+\frac{c_0\tau}2-V(x)}S,
\]
and letting $S\to\infty$ yields condition (ii). Conditions (i) and (ii) imply \eqref{eq:af-bound} with $B=\frac\beta2-c_0$.
\end{proof}

\begin{remark}\label{rem:direct-proof}
The functional $G_\tau$ lends itself to modification. For example, other weights can take the place of $1-\frac t\tau$ in \eqref{eq:window-functional}, and then a weighted average of $Q(X_t)$ replaces the uniform average in \eqref{eq:window-condition}. A suitable choice of the weight may improve the resulting bound. We do not pursue this direction here.
\end{remark}

\section{Infinite Prandtl number convection}\label{sec:system}

\subsection{The system and the notion of a solution}\label{subsec:solutions}

\subsubsection{\texorpdfstring{Reduction to an equation for $\theta$}{Reduction to an equation for theta}}\label{subsubsec:system}

For Rayleigh--B\'enard convection, one typically makes the transform $\theta\coloneqq T-(1-z)$, with $\theta:D\times[0,\infty)\to\R$, for the fluctuations around the conductive profile. Making this substitution in the system \eqref{eq:T-system}--\eqref{eq:zero-mean}, we obtain
\begin{equation}\label{eq:theta-system}
\begin{aligned}
    &{-}\Delta\vu+\nabla p=\Ray\,\theta\ez,\\
    &\divg\vu=0,\\
    &d\theta=\bigl[\Delta\theta-\vu\cdot\nabla\theta+u_2\bigr]\,dt+\sigma\sum_n\zeta_n\,dW^n_t,
\end{aligned}
\end{equation}
with
\begin{equation}\label{eq:theta-bc}
\begin{aligned}
    &\theta=u_2=\partial_zu_1=0\quad\text{at }z=0,1,\\
    &\int_0^1\!\int_0^Lu_1\,dx\,dz=0.
\end{aligned}
\end{equation}
Here we replaced $p$ by $p-\Ray(z-z^2/2)$, which absorbs the gradient $\Ray(1-z)\ez=\nabla\bigl(\Ray(z-z^2/2)\bigr)$ into the pressure.

At infinite Prandtl number, the first equation of \eqref{eq:theta-system} has no time derivative, and therefore at each instant $\vu$ is completely determined by $\theta$. We recall the standard facts we use about $\vu$ as a function of $\theta$ (see \cite[Section~2]{DC01} for no-slip boundaries).

For $\theta\in L^2(D)$, the Stokes system formed by the first two equations of \eqref{eq:theta-system}, with the conditions on $\vu$ in \eqref{eq:theta-bc}, has, by the Lax--Milgram theorem, a unique weak solution $\vu[\theta]=(u_1[\theta],u_2[\theta])$, with $\vu[\theta]$ and $\nabla\vu[\theta]$ in $L^2(D)$, and $\vu[\theta]$ is linear in $\theta$. Since $\divg\vu[\theta]=0$, $\vu[\theta]=(-\partial_z\psi,\partial_x\psi)$ for a stream function $\psi:D\to\R$. Taking the curl of the first equation of \eqref{eq:theta-system} eliminates $p$ and yields
\begin{equation}\label{eq:stream}
\begin{aligned}
    &\Delta^2\psi=-\Ray\,\partial_x\theta,\\
    &\psi=\partial_z^2\psi=0\quad\text{at }z=0,1.
\end{aligned}
\end{equation}
Here $\psi$ can be chosen to vanish at both $z=0$ and $z=1$ because of the zero-mean condition in \eqref{eq:theta-bc}.
Testing the first equation of \eqref{eq:theta-system}, the equation for the velocity, with $\vu[\theta]$ yields
\begin{equation}\label{eq:stokes-energy}
    \norm{\nabla\vu[\theta]}_{L^2(D)}^2=\Ray\,(u_2[\theta],\theta).
\end{equation}

Since $\vu[\theta]$ is linear in $\theta$, we denote the operator $\theta\mapsto\frac1\Ray u_2[\theta]$ by $K:L^2(D)\to L^2(D)$:
\begin{equation}\label{eq:K-def}
    u_2=\Ray\,K\theta .
\end{equation}
Hence, \eqref{eq:theta-system}--\eqref{eq:theta-bc} becomes one equation in $\theta$:
\begin{equation}\label{eq:theta-equation}
\begin{aligned}
    &d\theta=\bigl[\Delta\theta-\vu[\theta]\cdot\nabla\theta+\Ray K\theta\bigr]\,dt+\sigma\sum_n\zeta_n\,dW^n_t,\\
    &\theta=0\quad\text{at }z=0,1.
\end{aligned}
\end{equation}

\subsubsection{\texorpdfstring{The eigenbasis of $\Delta$ and $K$}{The eigenbasis of Delta and K}}\label{subsubsec:diagonal}

Note that each basis vector in \eqref{eq:fourier-basis} is $0$ at $z=0,1$. For $\theta=c_{k,n}=\frac2{L^{1/2}}\cos(kx)\sin(n\pi z)$, $\partial_x\theta=-k\,s_{k,n}$ and $\Delta^2s_{k,n}=(k^2+n^2\pi^2)^2s_{k,n}$, and hence $\psi$ in \eqref{eq:stream}, with $\Delta^2\psi=-\Ray\,\partial_x\theta$, is
\[
    \psi=\frac{\Ray\,k}{(k^2+n^2\pi^2)^2}\,s_{k,n},\qquad u_2=\partial_x\psi=\frac{\Ray\,k^2}{(k^2+n^2\pi^2)^2}\,c_{k,n}.
\]
The case $\theta=s_{k,n}$ is the same, and for $\theta=e_n$, $\partial_x\theta=0$, and therefore $\psi=0$. Consequently, $\Delta$ and $K$ are diagonal in the basis $\{e_n,c_{k,n},s_{k,n}\}$:
\begin{equation}\label{eq:fourier-operators}
\begin{aligned}
    \Delta c_{k,n}&=-(k^2+n^2\pi^2)\,c_{k,n},&\qquad Kc_{k,n}&=\frac{k^2}{(k^2+n^2\pi^2)^2}\,c_{k,n},\\
    \Delta s_{k,n}&=-(k^2+n^2\pi^2)\,s_{k,n},&\qquad Ks_{k,n}&=\frac{k^2}{(k^2+n^2\pi^2)^2}\,s_{k,n},\\
    \Delta e_n&=-n^2\pi^2\,e_n,&\qquad Ke_n&=0 .
\end{aligned}
\end{equation}
Thus, $K=-\partial_x^2\Delta^{-2}$ in this basis, and the eigenvalues of $K$ are nonnegative.

We also need an estimate on $\norm\vu_{L^\infty(D)}$, for the term $\vu\cdot\nabla\theta$ in \eqref{eq:theta-equation}. Let $\eta\coloneqq\Delta\psi$. By \eqref{eq:stream},
\[
\begin{aligned}
    \Delta\eta&=-\Ray\,\partial_x\theta,\qquad&\eta&=0\quad\text{at }z=0,1,\\
    \Delta\psi&=\eta,&\psi&=0\quad\text{at }z=0,1.
\end{aligned}
\]
For $m\geq1$, let $H^m(D)$ be the Sobolev space of functions on $D$, $L$-periodic in $x$, with derivatives up to order $m$ in $L^2(D)$. Since $\vu=(-\partial_z\psi,\partial_x\psi)$, the Sobolev embedding $H^2(D)\subset L^\infty(D)$, elliptic regularity for these two Dirichlet problems, and the estimate $\norm{\partial_x\theta}_{H^{-1}(D)}\leq\norm\theta_{L^2(D)}$ from \eqref{eq:Hminus1-fourier} imply
\begin{equation}\label{eq:u-bound}
    \norm\vu_{L^\infty(D)}\leq C\norm\vu_{H^2(D)}\leq C\norm\psi_{H^3(D)}\leq C\norm\eta_{H^1(D)}\leq C\Ray\norm{\partial_x\theta}_{H^{-1}(D)}\leq C\Ray\norm\theta_{L^2(D)}.
\end{equation}
By \eqref{eq:fourier-operators}, since $(k^2+n^2\pi^2)^2\geq4k^2n^2\pi^2$, the eigenvalues of $K$ are $0$ and
\begin{equation}\label{eq:K-eig-bound}
    \frac{k^2}{(k^2+n^2\pi^2)^2}\leq\frac{k^2}{4k^2n^2\pi^2}\leq\frac1{4\pi^2},
\end{equation}
and therefore, for $\theta\in L^2(D)$,
\begin{equation}\label{eq:u2-theta-bound}
    (u_2,\theta)=\Ray(K\theta,\theta)\leq\frac{\Ray}{4\pi^2}\norm{\theta}_{L^2(D)}^2.
\end{equation}

\subsubsection{The notion of a solution}\label{subsubsec:solution}

Lemma~\ref{lem:window} is stated for a Markov process, and consequently we need solutions of \eqref{eq:theta-equation} that exist, are unique and form a Markov process on $L^2(D)$.

We work on a filtered probability space $(\Omega,\mathcal F,\mathcal F_{t\geq0},\mathbb P)$ that is complete and satisfies the usual conditions, on which the $W^n$ are independent standard Brownian motions with respect to $\mathcal F_t$.

Let $H^1_0(D)$ and $H^{-1}(D)$ be as in \eqref{eq:H1-fourier} and \eqref{eq:Hminus1-fourier}, and let $\langle\cdot,\cdot\rangle_{H^{-1}\times H^1_0}$ be the pairing between $H^{-1}(D)$ and $H^1_0(D)$. For a general $\theta\in H^1_0(D)$, the term $\Delta\theta$ in \eqref{eq:theta-equation} is defined in the sense of distributions, as an element of $H^{-1}(D)$, by $\langle\Delta\theta,\vartheta\rangle_{H^{-1}\times H^1_0}=-(\nabla\theta,\nabla\vartheta)$ for $\vartheta\in H^1_0(D)$. Also, since $\divg\vu[\theta]=0$, integration by parts yields $(\vu[\theta]\cdot\nabla\theta,\vartheta)=-(\theta,\vu[\theta]\cdot\nabla\vartheta)$ for $\theta,\vartheta\in H^1_0(D)$. Testing \eqref{eq:theta-equation} against $\vartheta\in H^1_0(D)$ and integrating by parts leads to the following notion of a weak solution.

\begin{definition}\label{def:solution}
Let $\theta_0\in L^2(D)$ be deterministic. A \emph{solution} of \eqref{eq:theta-equation} with $\theta(0)=\theta_0$ is an $\mathcal F_t$-adapted process $\theta:[0,\infty)\times\Omega\to L^2(D)$ such that, almost surely, $\theta\in C([0,t];L^2(D))\cap L^2(0,t;H^1_0(D))$ for every $t>0$, and, for all $t\geq0$ and $\vartheta\in H^1_0(D)$,
\begin{equation}\label{eq:weak-solution}
\begin{aligned}
    (\theta(t),\vartheta)={}&(\theta_0,\vartheta)+\int_0^t\Bigl(-(\nabla\theta(r),\nabla\vartheta)+(\theta(r),\vu[\theta(r)]\cdot\nabla\vartheta)+\Ray(K\theta(r),\vartheta)\Bigr)\,dr\\
    &+\sigma\sum_n(\zeta_n,\vartheta)W^n_t,
\end{aligned}
\end{equation}
and such that, for every $t>0$,
\begin{equation}\label{eq:solution-moments}
    \E_{\theta_0}\int_0^t\norm{\nabla\theta(r)}_{L^2(D)}^2\,dr<\infty.
\end{equation}
\end{definition}

\begin{proposition}[Well-posedness and the Markov property]\label{prop:wellposed}
For deterministic $\theta_0\in L^2(D)$, \eqref{eq:theta-equation} has a unique solution $\theta$ with $\theta(0)=\theta_0$, in the sense of Definition~\ref{def:solution}. Moreover:
\begin{enumerate}[label=(\roman*)]
    \item The solutions form a time-homogeneous Markov process on $L^2(D)$, and, for every $t>0$,
    \[
        \E_{\theta_0}\int_0^t\Bigl(\norm{\nabla\theta(r)}_{L^2(D)}^2+\norm{\theta(r)}_{L^2(D)}^4\Bigr)\,dr<\infty.
    \]
    In particular, $\E_{\theta_0}\int_0^t\norm{\theta(r)}_{L^2(D)}^2\,dr<\infty$.
    \item The transition kernel of this Markov process is jointly measurable in $(t,\theta_0)$.
    \item It\^o's formula for the squared norm holds: for every $t\geq0$, almost surely,
    \begin{equation}\label{eq:energy-ito}
    \begin{aligned}
        \norm{\theta(t)}_{L^2(D)}^2&=\norm{\theta_0}_{L^2(D)}^2+\int_0^t\Bigl(-2\norm{\nabla\theta(r)}_{L^2(D)}^2+2(u_2[\theta(r)],\theta(r))+\sigma^2\Bigr)\,dr\\
        &\quad+2\sigma\sum_n\int_0^t(\theta(r),\zeta_n)\,dW^n_r ,
    \end{aligned}
    \end{equation}
    and the stochastic integral is a martingale.
\end{enumerate}
\end{proposition}

The proof is in Appendix~\ref{sec:wellposed}.

\subsection{The auxiliary functional}\label{subsec:aux-functional}

\subsubsection{The heat transport}\label{subsubsec:transport}

Since $T=\theta+1-z$, at time $t$,
\begin{equation}\label{eq:transport-F}
    \frac1L\int_0^1\!\!\int_0^L(u_2T-\partial_zT)\,dx\,dz=\frac1L\int_0^1\!\!\int_0^L\bigl(u_2\theta+u_2(1-z)+1-\partial_z\theta\bigr)\,dx\,dz=1+F(\theta(t)),
\end{equation}
because $\int_0^Lu_2\,dx=\int_0^L\partial_x\psi\,dx=0$ for every $z\in[0,1]$, by periodicity in $x$, and $\int_0^1\partial_z\theta\,dz=0$, as $\theta=0$ at $z=0,1$. Here $F:L^2(D)\to\R$ is
\begin{equation}\label{eq:observable}
    F(\theta)\coloneqq\frac1L(u_2,\theta)=\frac1L\int_0^1\!\int_0^L u_2(x,z)\theta(x,z)\,dx\,dz .
\end{equation}
Therefore, for the solution $\theta$ with $\theta(0)=\theta_0=T(0)-(1-z)$, \eqref{eq:nusselt} reads
\begin{equation}\label{eq:lower-time-average}
    \underline{\Nuss}-1=\liminf_{t\to\infty}\frac1t\,\E_{\theta_0}\int_0^tF(\theta(s))\,ds.
\end{equation}

We apply Lemma~\ref{lem:window} to \eqref{eq:theta-equation}, with $X_t\coloneqq\theta(t)$, with values in $\mathcal S\coloneqq L^2(D)$, and $F$ as in \eqref{eq:observable}.

\subsubsection{\texorpdfstring{The energy functional $V_0$}{The energy functional V0}}\label{subsubsec:generator}

Since Lemma~\ref{lem:window} needs a functional $V_0$, bounded above, such that \eqref{eq:balance} holds with some $Q\geq0$ and $c_0\geq0$, we choose a negative multiple of $\norm\theta_{L^2(D)}^2$,
\begin{equation}\label{eq:V0-def}
    V_0:L^2(D)\to\R,\qquad V_0(\theta)\coloneqq-\frac{\norm\theta_{L^2(D)}^2}{4L}\leq0.
\end{equation}
Then It\^o's formula \eqref{eq:energy-ito} for $\norm{\theta}_{L^2(D)}^2$ implies, in the form \eqref{eq:ito-drift},
\begin{equation}\label{eq:V0-ito}
    V_0(\theta(t))=V_0(\theta_0)+\int_0^t\LL V_0(\theta(r))\,dr+M_t,\qquad M_t\coloneqq-\frac\sigma{2L}\sum_n\int_0^t(\theta(r),\zeta_n)\,dW^n_r,
\end{equation}
with the drift
\begin{equation}\label{eq:LV0-def}
    \LL V_0(\theta)\coloneqq\frac1{2L}\norm{\nabla\theta}_{L^2(D)}^2-\frac1{2L}(u_2,\theta)-\frac{\sigma^2}{4L}\qquad\text{for }\theta\in H^1_0(D).
\end{equation}
The factor $\frac1{4L}$ is chosen so that adding $F=\frac1L(u_2,\theta)$ from \eqref{eq:observable} leaves $(u_2,\theta)$ with a positive coefficient: on $H^1_0(D)$,
\[
F+\LL V_0=\frac1{2L}(u_2,\theta)+\frac1{2L}\norm{\nabla\theta}_{L^2(D)}^2-\frac{\sigma^2}{4L}.
\]
The right side is $+\infty$ off $H^1_0(D)$, while Lemma~\ref{lem:window} needs real-valued functions. Thus, given $\beta>0$, we define
\begin{equation}\label{eq:Q-def}
    Q(\theta)\coloneqq\begin{cases}\dfrac1{2L}(u_2,\theta)+\dfrac1{2L}\norm{\nabla\theta}_{L^2(D)}^2&\text{for }\theta\in H^1_0(D),\\[1ex]\beta&\text{for }\theta\notin H^1_0(D),\end{cases}
\end{equation}
and we set
\begin{equation}\label{eq:LV0-off}
    \LL V_0(\theta)\coloneqq\beta-\frac{\sigma^2}{4L}-F(\theta)\qquad\text{for }\theta\notin H^1_0(D).
\end{equation}
Then $Q:L^2(D)\to[0,\infty)$, and \eqref{eq:balance} holds on $L^2(D)$ with $c_0\coloneqq\sigma^2/(4L)$. Here $Q\geq0$ because $(u_2,\theta)=\norm{\nabla\vu}_{L^2(D)}^2/\Ray\geq0$ by \eqref{eq:stokes-energy}.

The second case in \eqref{eq:Q-def} is only a technical convention. Almost surely, the solution satisfies $\theta(t)\in H^1_0(D)$ for almost every $t$, by Proposition~\ref{prop:wellposed}(i), and therefore the values of $Q$ and $\LL V_0$ off $H^1_0(D)$ change neither the formula \eqref{eq:V0-ito} for $V_0(\theta(t))$ nor any time integral of $Q(\theta(t))$. The value $\beta$ is chosen so that every $\theta$ with $Q(\theta)<\beta$ belongs to $H^1_0(D)$.

With $V_0$ alone, conditions (i) and (ii) in \eqref{eq:af-conditions} yield only $\underline{\Nuss}-1\geq-\sigma^2/(4L)$. Note that Lemma~\ref{lem:window} reduces Theorem~\ref{thm:lower-bound} to the proof of \eqref{eq:window-condition} for $Q$ over a fixed time $\tau$.

\begin{proposition}\label{prop:reduction}
Let $\beta,\tau>0$, and suppose that, for every deterministic $\theta_0\in H^1_0(D)$ with $Q(\theta_0)<\beta$, the solution $\theta$ of \eqref{eq:theta-equation} with $\theta(0)=\theta_0$ satisfies
\begin{equation}\label{eq:short-time-statement}
    \E_{\theta_0}\frac1\tau\int_0^\tau Q(\theta(t))\,dt\geq\beta .
\end{equation}
Then $\underline{\Nuss}-1\geq\frac\beta2-\frac{\sigma^2}{4L}$ for every deterministic initial temperature $T(0)$ with $T(0)-(1-z)\in L^2(D)$.
\end{proposition}

\begin{proof}
We wish to apply Lemma~\ref{lem:window} with $X_t=\theta(t)$, $\mathcal S=L^2(D)$, $F$ as in \eqref{eq:observable}, $V_0$ as in \eqref{eq:V0-def}, $\LL V_0$ as in \eqref{eq:LV0-def} and \eqref{eq:LV0-off}, $Q$ as in \eqref{eq:Q-def}, and $c_0=\sigma^2/(4L)$. For this, we check each assumption of Lemma~\ref{lem:window} in turn.

\noindent\textbf{The process $X_t$ is a time-homogeneous Markov process whose transition kernel is jointly measurable, and $(t,\omega)\mapsto X_t(\omega)$ is jointly measurable.}\\
For $X_t=\theta(t)$, the first two properties are Proposition~\ref{prop:wellposed}(i) and (ii), and the third holds because $\theta$ is adapted and has continuous paths in $L^2(D)$, by Definition~\ref{def:solution}.

\noindent\textbf{The function $F$ is Borel, $F\geq0$, and \eqref{eq:F-integrable} holds for every $x\in\mathcal S$.}\\
The function $F$ is continuous on $L^2(D)$ by \eqref{eq:u2-theta-bound}, $F=\frac1L(u_2,\theta)\geq0$ by \eqref{eq:stokes-energy}, and \eqref{eq:F-integrable} holds for every $x=\theta_0\in L^2(D)$, by Proposition~\ref{prop:wellposed}(i) and \eqref{eq:u2-theta-bound}.

\noindent\textbf{The functional $V_0$ is bounded above and has an extended generator $\LL V_0$.}\\
We have $V_0\leq0$. The function $\LL V_0$ is Borel on $L^2(D)$, because $H^1_0(D)$ is a Borel subset of $L^2(D)$, $\norm{\nabla\theta}_{L^2(D)}^2$ is a Borel function of $\theta\in H^1_0(D)$, and $F$ is continuous. Finally, $\LL V_0$ is an extended generator of $V_0$, as in \eqref{eq:ito-drift}, since, almost surely, $\theta(t)\in H^1_0(D)$ for almost every $t$, $\LL V_0$ satisfies \eqref{eq:LV-integrable} by Proposition~\ref{prop:wellposed}(i) and \eqref{eq:u2-theta-bound}, and $M_t$ in \eqref{eq:V0-ito} is a martingale by Proposition~\ref{prop:wellposed}(iii).

\noindent\textbf{The function $Q$ is Borel with values in $[0,\infty)$, and \eqref{eq:balance} holds on $\mathcal S$.}\\
The function $Q$ is Borel for the same reason as $\LL V_0$, and \eqref{eq:balance} holds on $L^2(D)$ with $Q$ as in \eqref{eq:Q-def} and $c_0=\sigma^2/(4L)$.

\noindent\textbf{Condition \eqref{eq:window-condition} holds for some $\tau>0$ independent of $x\in\mathcal S_{Q<\beta}$.}\\
Since $\mathcal S_{Q<\beta}\subset H^1_0(D)$ by \eqref{eq:Q-def}, this is the assumption \eqref{eq:short-time-statement} of the proposition.

Consequently, Lemma~\ref{lem:window} applies, and its conclusion, with \eqref{eq:lower-time-average}, is $\underline{\Nuss}-1\geq\frac\beta2-\frac{\sigma^2}{4L}$.
\end{proof}

\subsubsection{\texorpdfstring{The horizontal mean and the set $\mathcal S_{Q<\beta}$}{The horizontal mean and the set Omega}}\label{subsubsec:low-set}

We wish to show the estimate \eqref{eq:low-set-bounds} for $v\in\mathcal S_{Q<\beta}$. Let $v\in H^1_0(D)$, and recall that $\mathcal S_{Q<\beta}\subset H^1_0(D)$ by \eqref{eq:Q-def}. We split $v$ into its horizontal mean $\bar v:[0,1]\to\R$ and the rest:
\begin{equation}\label{eq:mean-def}
    \bar v(z)\coloneqq\frac1L\int_0^Lv\,dx,\qquad v'\coloneqq v-\bar v .
\end{equation}
Since the horizontal mean of $v'$ is $0$ by \eqref{eq:mean-def}, Fubini's theorem implies
\[
    \int_0^1\!\!\int_0^L\nabla\bar v\cdot\nabla v'\,dx\,dz=\int_0^1\partial_z\bar v\,\partial_z\Bigl(\int_0^Lv'\,dx\Bigr)dz=0,
\]
and hence
\[
    \norm{\nabla v}_{L^2(D)}^2=L\int_0^1(\partial_z\bar v)^2dz+\norm{\nabla v'}_{L^2(D)}^2.
\]
Also, since $\bar v$ depends only on $z$, $\bar v$ belongs to the closed span of the $e_n$. Recall that $Ke_n=0$ by \eqref{eq:fourier-operators}, and therefore $K\bar v=0$ and $(u_2[v],v)=\Ray(Kv,v)=\Ray(Kv',v')$. Hence, by the definition \eqref{eq:Q-def} of $Q$,
\begin{equation}\label{eq:mean-split}
    Q(v)=\frac\Ray{2L}(Kv',v')+\frac1{2L}\norm{\nabla v'}_{L^2(D)}^2+\frac12\int_0^1(\partial_z\bar v)^2dz .
\end{equation}
All three terms are nonnegative, and therefore
\begin{equation}\label{eq:low-set-bounds}
    \int_0^1(\partial_z\bar v)^2dz<2\beta\qquad\text{for }v\in\mathcal S_{Q<\beta}.
\end{equation}

\section{\texorpdfstring{Proof of Theorem~\ref{thm:lower-bound}}{Proof of Theorem 1.1}}\label{sec:proof-main}

By Proposition~\ref{prop:reduction}, it suffices to show that, for $\beta=c(\sigma^3\Ray)^{7/192}$ and some $\tau>0$, every deterministic $\theta_0\in H^1_0(D)$ with $Q(\theta_0)<\beta$ satisfies \eqref{eq:short-time-statement}, where $\theta$ is the solution of \eqref{eq:theta-equation} with $\theta(0)=\theta_0$. A direct estimate of the average of $Q$ in \eqref{eq:short-time-statement} is difficult. Instead, note that $Q\geq\frac1{2L}(u_2,\theta)=\frac\Ray{2L}(K\theta,\theta)$ on $H^1_0(D)$, and $(K\theta,\theta)$ is the squared length of $\theta$ for the symmetric nonnegative form
\begin{equation}\label{eq:K-form}
    (v,w)_K\coloneqq(Kv,w),\qquad v,w\in L^2(D).
\end{equation}
Consequently, the projection of $\theta$ on any direction already estimates $Q$ from below: for $\Psi\in L^2(D)$ with $(\Psi,\Psi)_K=1$, the projection
\[
    Y(t)\coloneqq\bigl(\theta(t),K\Psi\bigr)=\bigl(\theta(t),\Psi\bigr)_K
\]
satisfies $Y(t)^2\leq(K\theta(t),\theta(t))\leq\frac{2L}\Ray\,Q(\theta(t))$ whenever $\theta(t)\in H^1_0(D)$, by the Cauchy--Schwarz inequality for $(\cdot,\cdot)_K$. We therefore project the equation for $\theta$ onto the single direction $\Psi$, and estimate the time average of $Q$ from below by the time average of $Y^2$. More generally, we look for a real-valued process $Y$ and a constant $\gamma>0$ such that, almost surely, $Y(t)^2\leq\gamma\,Q(\theta(t))$ for almost every $t\in[0,\tau]$. Note that for such $Y$, if
\begin{equation*}
    \mathbb P_{\theta_0}\Bigl(\frac1\tau\int_0^\tau Y(t)^2\,dt>2\gamma\beta\Bigr)\geq\frac12,
\end{equation*}
then
\[
    \E_{\theta_0}\frac1\tau\int_0^\tau Q(\theta(t))\,dt\geq\frac1\gamma\cdot2\gamma\beta\cdot\frac12=\beta.
\]
The following lemma estimates the probability that the time average of $Y^2$ over $[0,\tau]$ is at most a given level $a^2$ and, at the same time, the error term in the equation for $Y$ is small:

\begin{lemma}\label{lem:gaussian}
Fix $a,\tau,\Lambda,q>0$. Let $Y$ be a real-valued It\^o process with deterministic $Y(0)$ satisfying
\begin{equation}\label{eq:scalar-sde}
    dY=\bigl(\Lambda Y+\epsilon(t)\bigr)\,dt+\sqrt q\,dW_t\qquad\text{on }[0,\tau],
\end{equation}
where $W$ is a standard Brownian motion, $\epsilon$ is an adapted process, and $\int_0^\tau\abs{\epsilon(t)}\,dt<\infty$ almost surely. If
\begin{equation}\label{cond:window}
    \Lambda\tau\geq1,\qquad e^{\Lambda\tau}\geq\frac{128\,a^2\Lambda}q,
\end{equation}
holds, then the event
\begin{equation}\label{eq:Gamma-Y}
    \Gamma_Y\coloneqq\Bigl\{\frac1\tau\int_0^\tau Y(t)^2\,dt\leq a^2\quad\text{and}\quad\frac{\int_0^\tau\abs{\epsilon(t)}\,dt}{(q/\Lambda)^{1/2}}\leq\frac18\Bigr\}
\end{equation}
satisfies $\mathbb P(\Gamma_Y)\leq\frac12$.
\end{lemma}

The proof of Lemma~\ref{lem:gaussian} is in Section~\ref{sec:escape-lemma}. The next proposition provides a process $Y$ to which Lemma~\ref{lem:gaussian} applies.

\begin{proposition}\label{prop:construct-Y}
There are constants $c,C>0$, depending only on $L$, such that the following holds. Let $0<\sigma\leq1$, and let $\sigma$ and $\Ray$ satisfy \eqref{eq:parameter-regime} with this constant $C$, and let $\beta\coloneqq c\,(\sigma^3\Ray)^{7/192}$. Then $\beta\geq1$ and $\beta\geq1/L$, and there are $\tau>0$ and $\gamma>0$ such that the following holds for every deterministic $\theta_0\in H^1_0(D)$ with $Q(\theta_0)<\beta$, for $Q$ in \eqref{eq:Q-def}. Let $\theta$ be the solution of \eqref{eq:theta-equation} with $\theta(0)=\theta_0$, and let
\begin{equation}\label{eq:Gamma-Q}
    \Gamma_Q\coloneqq\Bigl\{\frac1\tau\int_0^\tau Q(\theta(t))\,dt\leq2\beta\Bigr\}.
\end{equation}
Then there are $\Lambda>0$ and $q>0$, a standard Brownian motion $W$, an adapted process $\epsilon$ with $\int_0^\tau\abs{\epsilon(t)}\,dt<\infty$ almost surely, and a real-valued It\^o process $Y$ with deterministic $Y(0)$, such that, under $\mathbb P_{\theta_0}$, $Y$ satisfies \eqref{eq:scalar-sde}, \eqref{cond:window} holds with $a\coloneqq(2\gamma\beta)^{1/2}$, and
\begin{equation}\label{cond:error}
    \frac{\int_0^\tau\abs{\epsilon(t)}\,dt}{(q/\Lambda)^{1/2}}\leq\frac18
    \qquad\text{almost surely on }\Gamma_Q,
\end{equation}
\begin{equation}\label{cond:coordinate}
    Y(t)^2\leq\gamma\,Q(\theta(t))\quad\text{almost surely, for almost every }t\in[0,\tau].
\end{equation}
\end{proposition}

The proof of Proposition~\ref{prop:construct-Y} is in Section~\ref{subsec:proof-construct-Y}. This proof uses the process $Y$ defined in \eqref{eq:coordinate} and its equation \eqref{eq:coordinate-dynamics}.

\begin{proof}[Proof of Theorem~\ref{thm:lower-bound}]
Let $c$ and $C$ be the constants of Proposition~\ref{prop:construct-Y}, let $\sigma$ and $\Ray$ satisfy \eqref{eq:parameter-regime} with this $C$, and let $\beta$, $\tau$ and $\gamma$ be as in Proposition~\ref{prop:construct-Y}. Fix a deterministic $\theta_0\in H^1_0(D)$ with $Q(\theta_0)<\beta$, and let $\theta$, $\Gamma_Q$, $\Lambda$, $q$, $W$, $\epsilon$ and $Y$ be as in Proposition~\ref{prop:construct-Y}. On $\Gamma_Q$, \eqref{cond:coordinate} implies
\[
    \frac1\tau\int_0^\tau Y(t)^2\,dt\leq\frac\gamma\tau\int_0^\tau Q(\theta(t))\,dt\leq2\gamma\beta=a^2.
\]
Therefore, by \eqref{cond:error}, $\Gamma_Q$ is contained, up to a null set, in the event $\Gamma_Y$ in \eqref{eq:Gamma-Y}. Since $Y$ satisfies \eqref{eq:scalar-sde} and \eqref{cond:window} holds, Lemma~\ref{lem:gaussian} applies under $\mathbb P_{\theta_0}$ and yields $\mathbb P_{\theta_0}(\Gamma_Q)\leq\mathbb P_{\theta_0}(\Gamma_Y)\leq\frac12$. On $\Gamma_Q^c$, we have $\frac1\tau\int_0^\tau Q(\theta(t))\,dt>2\beta$, and $Q\geq0$, and consequently
\[
    \E_{\theta_0}\frac1\tau\int_0^\tau Q(\theta(t))\,dt\geq2\beta\bigl(1-\mathbb P_{\theta_0}(\Gamma_Q)\bigr)\geq\beta,
\]
which is \eqref{eq:short-time-statement}. This holds for every deterministic $\theta_0\in H^1_0(D)$ with $Q(\theta_0)<\beta$, and therefore Proposition~\ref{prop:reduction} implies
\begin{equation}\label{eq:final-bound}
    \underline{\Nuss}-1\geq\frac\beta2-\frac{\sigma^2}{4L}\geq\frac\beta4=\frac c4(\sigma^3\Ray)^{7/192},
\end{equation}
since $\sigma\leq1$ and $\beta\geq1/L$.
\end{proof}

\section{\texorpdfstring{Constructing $Y$}{Constructing Y}}\label{sec:coordinate}

In this section, we construct the process $Y$ of Proposition~\ref{prop:construct-Y}, which the proof of Theorem~\ref{thm:lower-bound} in Section~\ref{sec:proof-main} uses. In Section~\ref{subsec:spectral-setting}, we rewrite \eqref{eq:theta-equation} with a linear operator $\mathcal A_{\bar\theta}$, show that, for fixed $t$, $\mathcal A_{\bar\theta}$ is self-adjoint for an inner product $(\cdot,\cdot)_K$, and estimate its largest eigenvalue from below. In Section~\ref{subsec:construct-y}, we define $Y$ as the projection of the solution $\theta$ onto an eigenfunction for this eigenvalue, and derive the equation \eqref{eq:scalar-sde} for $Y$. Throughout this section, we assume
\begin{equation}\label{eq:construct-assumptions}
    \beta\geq1,\qquad\Ray\geq C(1+2\beta)^6.
\end{equation}

\subsection{Spectral setting}\label{subsec:spectral-setting}

We linearize \eqref{eq:theta-equation} about the horizontal mean $\bar\theta$ in \eqref{eq:mean-def}, because $\vu\cdot\nabla\bar\theta=u_2\partial_z\bar\theta$ is linear in $\theta$ once $\bar\theta$ is fixed, while $\vu\cdot\nabla\theta'$ is quadratic in $\theta'$.

With $\bar\theta$ and $\theta'$ as in \eqref{eq:mean-def}, since $u_2=\Ray K\theta$, we have $-\vu\cdot\nabla\theta+\Ray K\theta=-\vu\cdot\nabla\theta'+\Ray\,(1-\partial_z\bar\theta)K\theta$, and hence \eqref{eq:theta-equation} reads
\begin{equation}\label{eq:linearization}
\begin{aligned}
    &d\theta=\Bigl[\mathcal A_{\bar\theta}\theta-\vu\cdot\nabla\theta'\Bigr]dt+\sigma\sum_n\zeta_n\,dW^n_t ,\\
    &\mathcal A_{\bar\theta}\coloneqq\Delta+\Ray\,(1-\partial_z\bar\theta)K,
\end{aligned}
\end{equation}
where $\vu=\vu[\theta]$. Note that $\mathcal A_{\bar\theta}$ depends on $t$, since $\bar\theta=\bar\theta(t)$. For $v,\vartheta\in H^1_0(D)$,
\begin{equation}\label{eq:A-weak}
    \langle\mathcal A_{\bar\theta}v,\vartheta\rangle_{H^{-1}\times H^1_0}\coloneqq-(\nabla v,\nabla\vartheta)+\Ray\bigl((1-\partial_z\bar\theta)Kv,\vartheta\bigr),
\end{equation}
and $\vu\cdot\nabla\theta'\in L^2(D)$, since $\norm\vu_{L^\infty(D)}\leq C\Ray\norm\theta_{L^2(D)}$ by \eqref{eq:u-bound}.

The operator $\mathcal A_{\bar\theta}$ is not symmetric for $(\cdot,\cdot)$, and therefore we first find a Hilbert space on which $\mathcal A_{\bar\theta}$ is symmetric.

Fix a profile $\rho:[0,1]\to\R$, independent of time, with $\partial_z\rho\in L^2(0,1)$, and let
\begin{equation}\label{eq:A-rho}
    \mathcal A_\rho\coloneqq\Delta+\Ray\,(1-\partial_z\rho)K.
\end{equation} With $k_0=2\pi/L$ as in \eqref{eq:noise}, let
\begin{equation}\label{eq:H-def}
    \mathcal H\coloneqq\bigl\{\cos(k_0x)\,g(z):\ g\in L^2(0,1)\bigr\}=\overline{\operatorname{span}}\{c_{k_0,n}:n\geq1\},
\end{equation}
where the closure is in $L^2(D)$. We work in $\mathcal H$ because each $\zeta_n$ in \eqref{eq:noise} is a multiple of $c_{k_0,n}$. The operator $\mathcal A_\rho$ maps $\{v\in\mathcal H\cap H^1_0(D):\Delta v\in L^2(D)\}$ into $\mathcal H$. Indeed, for $v=\sum_na_nc_{k_0,n}\in\mathcal H\cap H^1_0(D)$ with $\Delta v\in L^2(D)$, that is, with $\sum_n(k_0^2+n^2\pi^2)^2a_n^2<\infty$, $\Delta v$ and $Kv$ belong to $\mathcal H$ by \eqref{eq:fourier-operators}, and so does $(1-\partial_z\rho)Kv$, since $\partial_z\rho\in L^2(0,1)$ and $Kv=\frac1\Ray u_2[v]$ is bounded, by \eqref{eq:u-bound}.

By \eqref{eq:fourier-operators}, $K$ is symmetric for $(\cdot,\cdot)$. In fact, observe that
\begin{equation}\label{eq:K-on-H}
    K=k_0^2\Delta^{-2}\qquad\text{on }\mathcal H,
\end{equation}
where $\Delta^{-1}$ denotes the inverse of $\Delta$ with boundary condition $0$ at $z=0,1$. Since $\Delta^{-1}$ is symmetric for $(\cdot,\cdot)$ and injective, the form $(\cdot,\cdot)_K$ in \eqref{eq:K-form} satisfies
\begin{equation}\label{eq:K-inner}
    (v,w)_K=k_0^2(\Delta^{-1}v,\Delta^{-1}w),\qquad v,w\in\mathcal H,
\end{equation}
and is an inner product on $\mathcal H$.

For $v,w\in\mathcal H\cap H^1_0(D)$ with $\Delta v,\Delta w\in L^2(D)$, the symmetry of $K$ and \eqref{eq:K-on-H} yield
\begin{equation}\label{eq:A-form}
\begin{aligned}
    (\mathcal A_\rho v,w)_K&=(K\Delta v,w)+\Ray\bigl((1-\partial_z\rho)Kv,Kw\bigr)\\
    &=k_0^2(\Delta^{-1}v,w)+\Ray\int_0^1\!\!\int_0^L(1-\partial_z\rho)(Kv)(Kw)\,dx\,dz.
\end{aligned}
\end{equation}
The right side of \eqref{eq:A-form} is symmetric in $v$ and $w$, because $\Delta^{-1}$ is symmetric for $(\cdot,\cdot)$. Thus, $\mathcal A_\rho$ is symmetric for $(\cdot,\cdot)_K$. The space $\mathcal H$ is not complete for $(\cdot,\cdot)_K$, and the spectral theorem needs a complete space. For this reason, the following proposition is stated on the completion of $\mathcal H$.

\begin{proposition}\label{prop:spectral}
Let $\mathcal H_K$ be the completion of $\mathcal H$ for $(\cdot,\cdot)_K$, with inner product also denoted by $(\cdot,\cdot)_K$. For every $\rho:[0,1]\to\R$ with $\partial_z\rho\in L^2(0,1)$, the operator $\mathcal A_\rho$ on $\{v\in\mathcal H\cap H^1_0(D):\Delta v\in L^2(D)\}$ extends to a self-adjoint operator on $\mathcal H_K$, with domain $\mathcal H$ and compact resolvent. This extension, also denoted by $\mathcal A_\rho$, is bounded above.
\end{proposition}

The proof is in Section~\ref{subsec:proof-spectral}. By Proposition~\ref{prop:spectral} and the spectral theorem for self-adjoint operators with compact resolvent, $\mathcal A_\rho$ has a $(\cdot,\cdot)_K$-orthonormal basis of eigenfunctions, which lie in its domain $\mathcal H$. Every eigenfunction $v$ of $\mathcal A_\rho$, with eigenvalue $\lambda$, satisfies $\Delta v=\lambda v-\Ray(1-\partial_z\rho)Kv\in\mathcal H$, and therefore $\Delta v\in L^2(D)$. The largest eigenvalue of $\mathcal A_\rho$ is
\begin{equation}\label{eq:largest-eigenvalue}
    \sup_{v\in\mathcal H,\,v\neq0}\frac{(\mathcal A_\rho v,v)_K}{(v,v)_K}.
\end{equation}

\begin{proposition}\label{prop:growth} Let $\theta\in H^1_0(D)$ have horizontal mean $\rho=\frac1L\int_0^L\theta\,dx$ as in \eqref{eq:mean-def}, and let $\mathcal A_\rho=\Delta+\Ray(1-\partial_z\rho)K$ be as in \eqref{eq:A-rho}. Suppose $\int_0^1(\partial_z\rho)^2dz<2\beta$ (which holds for $\theta\in\mathcal S_{Q<\beta}$, by \eqref{eq:low-set-bounds}) and $\Ray\geq C(1+2\beta)^6$, as in \eqref{eq:construct-assumptions}. Then the largest eigenvalue of $\mathcal A_\rho$ is at least
\begin{equation}\label{eq:growth-bound}
    \lambda_*\coloneqq c\,\Ray\,(1+2\beta)^{-4}.
\end{equation}
Here $c>0$ is fixed, depending only on $L$.
\end{proposition}

The proof is in Section~\ref{subsec:proof-growth}.

\subsection{\texorpdfstring{The projection $Y$}{The projection Y}}\label{subsec:construct-y}

For a starting state $\theta_0\in\mathcal S_{Q<\beta}$, with horizontal mean $\bar\theta_0$, let $\Lambda$ be the largest eigenvalue of $\mathcal A_{\bar\theta_0}$. By \eqref{eq:low-set-bounds} and Proposition~\ref{prop:growth} with $\rho=\bar\theta_0$, we have $\Lambda\geq\lambda_*$. By Proposition~\ref{prop:spectral}, choose an eigenfunction $\Psi\in\mathcal H$ for $\Lambda$ with $(\Psi,\Psi)_K=1$, that is,
\begin{equation}\label{eq:Psi-def}
    \mathcal A_{\bar\theta_0}\Psi=\Lambda\Psi,\qquad(\Psi,\Psi)_K=1,
\end{equation}
and let $\theta$ be the solution of \eqref{eq:theta-equation} with $\theta(0)=\theta_0$. Set
\begin{equation}\label{eq:coordinate}
    Y:[0,\infty)\times\Omega\to\R,\qquad Y(t)\coloneqq\bigl(\theta(t),K\Psi\bigr)=\bigl(P_{\mathcal H}\theta(t),\Psi\bigr)_K.
\end{equation}
Here $K$ and $(v,w)_K=(Kv,w)$ are as in \eqref{eq:K-def} and \eqref{eq:K-form}, and
\begin{equation}\label{eq:PH-def}
    P_{\mathcal H}\theta\coloneqq\sum_n(\theta,c_{k_0,n})\,c_{k_0,n}
\end{equation}
is the $L^2(D)$-orthogonal projection of $\theta$ onto $\mathcal H$. The second equality in \eqref{eq:coordinate} holds because $(P_{\mathcal H}\theta,w)_K=(\theta,Kw)$ for $w\in\mathcal H$, since $Kw\in\mathcal H$.

The noise in the equation for $Y=(\theta,K\Psi)$ is $\sigma\sum_n(\zeta_n,K\Psi)\,dW^n_t$, with $\zeta_n$ as in \eqref{eq:noise}. By \eqref{eq:noise} and then \eqref{eq:fourier-operators}, the variance per unit time of this noise is
\begin{equation}\label{eq:q-formula}
\begin{aligned}
q&\coloneqq\sigma^2\sum_n(\zeta_n,K\Psi)^2
  =\frac L4C_L^2\sigma^2\sum_n\frac{k_0^2}{(k_0^2+n^2\pi^2)^2}\,(c_{k_0,n},K\Psi)^2\\
  &=\frac L4C_L^2\sigma^2\,\bigl(K(K\Psi),K\Psi\bigr)
  =\frac L4C_L^2\sigma^2\,(K\Psi,K\Psi)_K ,
\end{aligned}
\end{equation}
and hence $q>0$. The choices $\Psi$, $\Lambda$ and $q$ may depend on $\theta_0$, but remain fixed in time.

To derive the equation for $Y$, we rewrite \eqref{eq:linearization} as
\begin{equation}\label{eq:linearization-frozen}
    d\theta=\Bigl[\mathcal A_{\bar\theta_0}\theta+\bigl(\mathcal A_{\bar\theta(t)}-\mathcal A_{\bar\theta_0}\bigr)\theta-\vu\cdot\nabla\theta'\Bigr]dt+\sigma\sum_n\zeta_n\,dW^n_t,
\end{equation}
where the operator $\mathcal A_{\bar\theta_0}$ does not depend on $t$.

\begin{proposition}\label{prop:coordinate}
Let $\beta$ and $\Ray$ satisfy \eqref{eq:construct-assumptions}, let $\theta_0\in\mathcal S_{Q<\beta}$, let $\Lambda$ be the largest eigenvalue of $\mathcal A_{\bar\theta_0}$, let $\Psi$ be as in \eqref{eq:Psi-def}, let $\theta$ be the solution of \eqref{eq:theta-equation} with $\theta(0)=\theta_0$, and let $Y$ be as in \eqref{eq:coordinate} and $q$ as in \eqref{eq:q-formula}. Then $Y(t)^2\leq\frac{2L}\Ray\,Q(\theta(t))$ for all $t\geq0$ with $\theta(t)\in H^1_0(D)$, and therefore, almost surely, for almost every $t$. Moreover,
\begin{equation}\label{eq:W-def}
    W\coloneqq\sigma q^{-1/2}\sum_n(\zeta_n,K\Psi)W^n
\end{equation}
is a standard Brownian motion, and
\begin{equation}\label{eq:coordinate-dynamics}
    dY=\bigl(\Lambda Y+\epsilon_1(t)+\epsilon_2(t)\bigr)\,dt+\sqrt q\,dW_t,
\end{equation}
that is, almost surely,
\[
    Y(t)=Y(0)+\int_0^t(\Lambda Y+\epsilon_1+\epsilon_2)\,dr+\sqrt q\,W_t\qquad\text{for all }t\geq0,
\]
where, for almost every $t$,
\begin{equation}\label{eq:eps-def}
\begin{aligned}
    \epsilon_1(t)&\coloneqq\Ray\bigl((\partial_z\bar\theta_0-\partial_z\bar\theta(t))KP_{\mathcal H}\theta(t),K\Psi\bigr),\\
    \epsilon_2(t)&\coloneqq\bigl(\theta'(t),\vu(t)\cdot\nabla K\Psi\bigr).
\end{aligned}
\end{equation}
Here $\bar\theta=\frac1L\int_0^L\theta\,dx$ and $\theta'=\theta-\bar\theta$ are as in \eqref{eq:mean-def}.
\end{proposition}

\begin{proof}[Proof of Proposition~\ref{prop:coordinate}]
\textbf{The estimate $Y^2\leq\frac{2L}\Ray Q$.} Fix $t\geq0$ with $\theta(t)\in H^1_0(D)$, and denote $\theta=\theta(t)$ and $Y=Y(t)$. Denote $P_{\mathcal H}^\perp\theta\coloneqq\theta-P_{\mathcal H}\theta$, which is orthogonal to $\mathcal H$ in $L^2(D)$. By \eqref{eq:fourier-operators}, $K$ maps $\mathcal H$ into $\mathcal H$, and the eigenvalues of $K$ are nonnegative. Therefore, $(P_{\mathcal H}\theta,P_{\mathcal H}^\perp\theta)_K=(KP_{\mathcal H}\theta,P_{\mathcal H}^\perp\theta)=0$ and $(P_{\mathcal H}^\perp\theta,P_{\mathcal H}^\perp\theta)_K\geq0$, and consequently
\begin{equation}\label{eq:PH-K-bound}
    (K\theta,\theta)=(\theta,\theta)_K=(P_{\mathcal H}\theta,P_{\mathcal H}\theta)_K+(P_{\mathcal H}^\perp\theta,P_{\mathcal H}^\perp\theta)_K\geq(P_{\mathcal H}\theta,P_{\mathcal H}\theta)_K.
\end{equation}
The estimate \eqref{eq:PH-K-bound}, the Cauchy--Schwarz inequality for $(\cdot,\cdot)_K$, $u_2=\Ray K\theta$ and $Q\geq\frac1{2L}(u_2,\theta)$, a consequence of \eqref{eq:Q-def}, imply
\begin{equation}\label{eq:coordinate-bound}
    Y^2\leq(P_{\mathcal H}\theta,P_{\mathcal H}\theta)_K\,(\Psi,\Psi)_K\leq(K\theta,\theta)=\frac{(u_2,\theta)}\Ray\leq\frac{2L}\Ray\,Q .
\end{equation}

\textbf{The Brownian motion $W$.} Denote $b_n\coloneqq\sigma q^{-1/2}(\zeta_n,K\Psi)$. Then $W=\sum_nb_nW^n$ by \eqref{eq:W-def}, and $\sum_nb_n^2=1$ by the definition of $q$ in \eqref{eq:q-formula}. Since the $W^n$ are independent standard Brownian motions, each partial sum $\sum_{n\leq N}b_nW^n$ has continuous paths and independent, centered Gaussian increments, with variance $(t-s)\sum_{n\leq N}b_n^2$ over a time interval $[s,t]$. To pass these properties to $W$, we show that the partial sums converge uniformly on bounded time intervals. For $N_1<N_2$, $t_1>0$ and $\lambda>0$, Doob's inequality, applied to the real-valued martingale $\sum_{N_1<n\leq N_2}b_nW^n$, yields
\[
    \mathbb P_{\theta_0}\Bigl(\sup_{t\leq t_1}\Bigl|\sum_{N_1<n\leq N_2}b_nW^n_t\Bigr|\geq\lambda\Bigr)\leq\frac1{\lambda^2}\,\E_{\theta_0}\Bigl(\sum_{N_1<n\leq N_2}b_nW^n_{t_1}\Bigr)^2=\frac{t_1}{\lambda^2}\sum_{N_1<n\leq N_2}b_n^2,
\]
which tends to $0$ as $N_1,N_2\to\infty$. Hence, the series defining $W$ converges in probability, uniformly on $[0,t_1]$. Since $t_1>0$ is arbitrary, the limit $W$ has continuous paths and independent, centered Gaussian increments, with variance $(t-s)\sum_nb_n^2=t-s$ over a time interval $[s,t]$, and therefore $W$ is a standard Brownian motion.

\textbf{The equation \eqref{eq:coordinate-dynamics} for $Y$.} By It\^o's formula for the linear functional $\theta\mapsto(\theta,K\Psi)$, that is, by \eqref{eq:weak-solution} with the fixed test function $\vartheta=K\Psi\in H^1_0(D)$, the equation \eqref{eq:linearization-frozen} yields
\begin{equation}\label{eq:dY-raw}
\begin{aligned}
    dY&=\Bigl[\langle\mathcal A_{\bar\theta_0}\theta,K\Psi\rangle_{H^{-1}\times H^1_0}+\bigl((\mathcal A_{\bar\theta(t)}-\mathcal A_{\bar\theta_0})\theta,K\Psi\bigr)-(\vu\cdot\nabla\theta',K\Psi)\Bigr]\,dt\\
    &\qquad+\sigma\sum_n(\zeta_n,K\Psi)\,dW^n_t,
\end{aligned}
\end{equation}
where, by \eqref{eq:A-weak} with $\bar\theta=\bar\theta_0$, $v=\theta$ and $\vartheta=K\Psi$,
\[
    \langle\mathcal A_{\bar\theta_0}\theta,K\Psi\rangle_{H^{-1}\times H^1_0}=-(\nabla\theta,\nabla K\Psi)+\Ray\bigl((1-\partial_z\bar\theta_0)K\theta,K\Psi\bigr).
\]
Recall that $\Delta\Psi\in L^2(D)$, because $\Psi$ is an eigenfunction of $\mathcal A_{\bar\theta_0}$. Integration by parts, the identity $\Delta K\Psi=K\Delta\Psi$ from \eqref{eq:fourier-operators}, the symmetry of $K$ for $(\cdot,\cdot)$, the definition of $\mathcal A_{\bar\theta_0}$ in \eqref{eq:linearization}, and $\mathcal A_{\bar\theta_0}\Psi=\Lambda\Psi$ from \eqref{eq:Psi-def} imply
\begin{equation}\label{eq:drift-linear}
\begin{aligned}
    \langle\mathcal A_{\bar\theta_0}\theta,K\Psi\rangle_{H^{-1}\times H^1_0}&=(\theta,K\Delta\Psi)+\Ray\bigl(\theta,K[(1-\partial_z\bar\theta_0)K\Psi]\bigr)\\
    &=(\theta,K\mathcal A_{\bar\theta_0}\Psi)=\Lambda(\theta,K\Psi)=\Lambda Y.
\end{aligned}
\end{equation}
By the definition of $\mathcal A_{\bar\theta}$ in \eqref{eq:linearization},
\begin{equation}\label{eq:drift-difference}
\begin{aligned}
    \bigl((\mathcal A_{\bar\theta(t)}-\mathcal A_{\bar\theta_0})\theta,K\Psi\bigr)&=\Ray\bigl((\partial_z\bar\theta_0-\partial_z\bar\theta(t))K\theta,K\Psi\bigr)\\
    &=\Ray\bigl((\partial_z\bar\theta_0-\partial_z\bar\theta(t))KP_{\mathcal H}\theta,K\Psi\bigr)=\epsilon_1(t),
\end{aligned}
\end{equation}
where the second equality holds because $(\partial_z\bar\theta_0-\partial_z\bar\theta(t))K\Psi$ belongs to $\mathcal H$ and $K(\theta-P_{\mathcal H}\theta)$ is orthogonal to $\mathcal H$ in $L^2(D)$, by \eqref{eq:fourier-operators}. Since $\divg\vu=0$, $\vu\cdot\nabla\theta'=\divg(\vu\theta')$, and integration by parts yields
\begin{equation}\label{eq:drift-nonlinear}
\begin{aligned}
    -(\vu\cdot\nabla\theta',K\Psi)&=-\int_0^1\!\!\int_0^L\divg(\vu\theta')K\Psi\,dx\,dz\\
    &=\int_0^1\!\!\int_0^L\theta'\,\vu\cdot\nabla K\Psi\,dx\,dz=(\theta',\vu\cdot\nabla K\Psi)=\epsilon_2(t).
\end{aligned}
\end{equation}
Substituting \eqref{eq:drift-linear}, \eqref{eq:drift-difference} and \eqref{eq:drift-nonlinear} into \eqref{eq:dY-raw}, and using $\sigma\sum_n(\zeta_n,K\Psi)\,dW^n_t=\sqrt q\,dW_t$, which holds by the definition \eqref{eq:W-def} of $W$, leads to \eqref{eq:coordinate-dynamics}.
\end{proof}

By Proposition~\ref{prop:coordinate}, $Y$ satisfies \eqref{cond:coordinate}, $Y(t)^2\leq\gamma\,Q(\theta(t))$, with
\begin{equation}\label{eq:gamma-def}
    \gamma\coloneqq\frac{2L}\Ray.
\end{equation}
Both \eqref{cond:window} and \eqref{cond:error} contain the ratio $q/\Lambda$, which the next proposition estimates from below.

\begin{proposition}\label{prop:noise}
With $\beta$ and $\Ray$ as in \eqref{eq:construct-assumptions}, $\theta_0\in\mathcal S_{Q<\beta}$, $\Lambda$ the largest eigenvalue of $\mathcal A_{\bar\theta_0}$, $\mathcal A_{\bar\theta_0}\Psi=\Lambda\Psi$ and $(\Psi,\Psi)_K=1$ as in \eqref{eq:Psi-def}, and $q=\sigma^2\sum_n(\zeta_n,K\Psi)^2$ as in \eqref{eq:q-formula},
\begin{equation}\label{eq:noise-bound}
\frac q\Lambda\geq\frac{c\,\sigma^2}\Ray\,(1+2\beta)^{-44/7}.
\end{equation}
\end{proposition}

The proof is provided in Section~\ref{subsec:proof-noise}.

\section{\texorpdfstring{Proof of Proposition~\ref{prop:construct-Y}}{Proof of Proposition 4.2}}\label{sec:verify}\label{subsec:proof-construct-Y}\label{subsec:choose-tau}\label{subsec:choose-beta}\label{sec:errors}

In this section, we show that the process $Y$ defined in \eqref{eq:coordinate} satisfies the conditions \eqref{cond:window} and \eqref{cond:error} of Proposition~\ref{prop:construct-Y}. This completes the proof of Proposition~\ref{prop:construct-Y}, and therefore the proof of Theorem~\ref{thm:lower-bound} in Section~\ref{sec:proof-main}, because the other two requirements of Proposition~\ref{prop:construct-Y}, the equation \eqref{eq:scalar-sde} and the estimate \eqref{cond:coordinate}, follow from Proposition~\ref{prop:coordinate}. We choose $\tau$ so that \eqref{cond:window} holds, we reduce \eqref{cond:error} to an estimate on the processes $\epsilon_1$ and $\epsilon_2$ of Proposition~\ref{prop:coordinate}, and we verify the conditions on $\beta$ that these two steps use. Throughout this section, $0<\sigma\leq1$.

We first consider any $\beta$ and $\Ray$ that satisfy \eqref{eq:construct-assumptions}, and return to $\beta=c\,(\sigma^3\Ray)^{7/192}$ in \eqref{eq:beta-choice}.

By \eqref{eq:low-set-bounds} and Proposition~\ref{prop:growth}, $\Lambda\geq\lambda_*$ for every $\theta_0\in\mathcal S_{Q<\beta}$, and therefore we choose $\tau$ independently of $\theta_0$ by
\begin{equation}\label{eq:window-length}
    \tau\coloneqq\frac\beta{\lambda_*}=\frac{\beta\,(1+2\beta)^4}{c\,\Ray},
\end{equation}
with the constant $c$ in \eqref{eq:growth-bound}. Since $\beta\geq1$ by \eqref{eq:construct-assumptions}, we have $\Lambda\tau\geq\lambda_*\tau=\beta\geq1$, which implies the first inequality in \eqref{cond:window}.

With $\gamma=\frac{2L}\Ray$ as in \eqref{eq:gamma-def} and $a^2=2\gamma\beta$, the second inequality in \eqref{cond:window} is equivalent to $e^{\Lambda\tau}\geq8^2\cdot4\gamma\beta\,\Lambda/q$, and Proposition~\ref{prop:noise} yields
\[
8^2\cdot\frac{4\gamma\beta\Lambda}q=\frac{8^3L\beta}{\Ray\,q/\Lambda}\leq C\sigma^{-2}\beta\,(1+2\beta)^{44/7}.
\]
Since $\Lambda\tau\geq\beta$, it therefore suffices that, for this $C$,
\begin{equation}
C\sigma^{-2}\beta\,(1+2\beta)^{44/7}\leq e^\beta,\label{eq:tau-sufficient}
\end{equation}
which is true if $\sigma^{-2}\leq e^{\beta/2}$ and $C\beta(1+2\beta)^{44/7}\leq e^{\beta/2}$, that is, if $\beta\geq4\log(1/\sigma)$ and $\beta\geq C$, for a constant $C$ depending only on $L$. Consequently, \eqref{cond:window} holds for
\begin{equation}\label{eq:beta-lower}
    \beta\geq C\bigl(1+\log(1/\sigma)\bigr).
\end{equation}

Next we turn to \eqref{cond:error}, for which we use the following estimate.

\begin{proposition}\label{prop:error-ratio}
Let $\beta$ and $\Ray$ satisfy \eqref{eq:construct-assumptions}, let $\theta_0\in\mathcal S_{Q<\beta}$, that is, $\theta_0\in H^1_0(D)$ with $Q(\theta_0)<\beta$ for $Q$ in \eqref{eq:Q-def}, let $\theta$ be the solution of \eqref{eq:theta-equation} with $\theta(0)=\theta_0$, let $\Lambda$ be the largest eigenvalue \eqref{eq:largest-eigenvalue} of $\mathcal A_{\bar\theta_0}$, let $\Psi$ be as in \eqref{eq:Psi-def}, $q$ as in \eqref{eq:q-formula}, and $\epsilon_1$ and $\epsilon_2$ as in \eqref{eq:eps-def}, that is,
\[
\begin{aligned}
    \epsilon_1(t)&=\Ray\bigl((\partial_z\bar\theta_0-\partial_z\bar\theta(t))KP_{\mathcal H}\theta(t),K\Psi\bigr),\\
    \epsilon_2(t)&=\bigl(\theta'(t),\vu(t)\cdot\nabla K\Psi\bigr),
\end{aligned}
\]
with $P_{\mathcal H}$ as in \eqref{eq:PH-def}, and $\bar\theta$ and $\theta'$ as in \eqref{eq:mean-def}, and let $\tau$ be as in \eqref{eq:window-length}. Then, almost surely,
\begin{equation}\label{eq:eps-integrable}
    \int_0^\tau\bigl(\abs{\epsilon_1(t)}+\abs{\epsilon_2(t)}\bigr)\,dt<\infty.
\end{equation}
Moreover, almost surely on the event $\Gamma_Q$ in \eqref{eq:Gamma-Q},
\begin{equation}\label{eq:error-ratio}
    \frac{\int_0^\tau\abs{\epsilon_1}\,dt+\int_0^\tau\abs{\epsilon_2}\,dt}{\bigl(q/\Lambda\bigr)^{1/2}}
    \leq\frac C\sigma\Bigl(\beta^{205/14}\Ray^{-2/3}+\beta^{64/7}\Ray^{-1/3}\Bigr).
\end{equation}
\end{proposition}

The proof is in Section~\ref{subsec:proof-errors}.
By Proposition~\ref{prop:error-ratio}, \eqref{cond:error} with $\epsilon=\epsilon_1+\epsilon_2$ holds if
\begin{equation}\label{eq:beta-conditions}
    C\sigma^{-1}\beta^{64/7}\Ray^{-1/3}\leq\frac1{16}\qquad\text{and}\qquad C\sigma^{-1}\beta^{205/14}\Ray^{-2/3}\leq\frac1{16}.
\end{equation}
Solving each inequality in \eqref{eq:beta-conditions} for $\beta$, and taking the smaller of the two resulting constants as $c$, we see that \eqref{eq:beta-conditions} holds if
\begin{align}
    \beta&\leq c\,(\sigma^3\Ray)^{7/192},\label{eq:beta-c-cond}\\
    \beta&\leq c\,\sigma^{14/205}\Ray^{28/615}.\label{eq:beta-c-cond2}
\end{align}
Observe that, for $\Ray$ large enough, \eqref{eq:beta-c-cond} implies \eqref{eq:beta-c-cond2}. Indeed, for $\sigma\leq1$ and $\Ray\geq1$,
\begin{equation}\label{eq:exponent-compare}
    (\sigma^3\Ray)^{7/192}=\sigma^{7/64}\Ray^{7/192}\leq\sigma^{14/205}\Ray^{28/615},
\end{equation}
because $\frac7{64}>\frac{14}{205}$ and $\frac7{192}<\frac{28}{615}$. Consequently, for $\sigma\leq1$ and $\Ray\geq1$, \eqref{eq:beta-conditions} holds if \eqref{eq:beta-c-cond} holds. In particular, $\beta=c\,(\sigma^3\Ray)^{7/192}$ satisfies \eqref{eq:beta-c-cond} for some sufficiently small $c$.

\begin{proof}[Proof of Proposition~\ref{prop:construct-Y}]
First, we fix the constants $c$ and $C$ of the proposition. Let $0<c\leq1$ be at most the constant in \eqref{eq:beta-c-cond}, and let
\begin{equation}\label{eq:beta-choice}
    \beta\coloneqq c\,(\sigma^3\Ray)^{7/192},
\end{equation}
so that \eqref{eq:beta-c-cond} holds. With this $c$ fixed, we choose the constant $C\geq1$ in \eqref{eq:parameter-regime} as follows. Since $\sigma\leq1$, \eqref{eq:parameter-regime} implies $\Ray\geq\sigma^3\Ray\geq C\geq1$, and therefore \eqref{eq:beta-c-cond} implies \eqref{eq:beta-conditions}. Moreover, \eqref{eq:parameter-regime} and \eqref{eq:beta-choice} imply
\[
    \beta\geq c\,C^{7/192}\bigl(1+\log(1/\sigma)\bigr),
\]
and therefore $\beta\geq1$, $\beta\geq1/L$ and \eqref{eq:beta-lower} hold if $C$ is sufficiently large. Finally, $\beta\leq\Ray^{7/192}$ by \eqref{eq:beta-choice}, and $\Ray\geq C$, and therefore \eqref{eq:construct-assumptions} also holds if $C$ is sufficiently large. We fix such a $C$. Consequently, $\beta$ in \eqref{eq:beta-choice} satisfies \eqref{eq:construct-assumptions}, \eqref{eq:beta-lower} and \eqref{eq:beta-conditions}, which are all the conditions on $\beta$ imposed so far.

Let $\tau\coloneqq\beta/\lambda_*$ as in \eqref{eq:window-length}, and let $\gamma=\frac{2L}\Ray$ as in \eqref{eq:gamma-def}. Fix a deterministic $\theta_0\in H^1_0(D)$ with $Q(\theta_0)<\beta$, that is, $\theta_0\in\mathcal S_{Q<\beta}$. Let $\Lambda$ be the largest eigenvalue \eqref{eq:largest-eigenvalue} of $\mathcal A_{\bar\theta_0}$, let $\Psi\in\mathcal H$ be an eigenfunction for $\Lambda$ with $(\Psi,\Psi)_K=1$, let $Y$ be as in \eqref{eq:coordinate} and $q$ as in \eqref{eq:q-formula}, let $W$, $\epsilon_1$ and $\epsilon_2$ be as in Proposition~\ref{prop:coordinate}, and let $\epsilon\coloneqq\epsilon_1+\epsilon_2$. Then $\Lambda\geq\lambda_*>0$ and $q>0$, and $Y(0)=(\theta_0,K\Psi)$ is deterministic. By Proposition~\ref{prop:coordinate}, $W$ is a standard Brownian motion, $Y$ satisfies \eqref{eq:scalar-sde}, and \eqref{cond:coordinate} holds.

The process $\epsilon$ is adapted, because $\epsilon_1(t)$ and $\epsilon_2(t)$ are functions of $\theta(t)$. By \eqref{eq:eps-integrable} in Proposition~\ref{prop:error-ratio},
\[
    \int_0^\tau\abs{\epsilon(t)}\,dt<\infty\qquad\text{almost surely.}
\]

Finally, \eqref{cond:window} with $a^2=2\gamma\beta$ holds because \eqref{eq:beta-lower} holds, and \eqref{cond:error} holds by Proposition~\ref{prop:error-ratio} and \eqref{eq:beta-conditions}.
\end{proof}

\section{Proofs of the technical results}\label{sec:proofs}

In this section, we prove Lemma~\ref{lem:gaussian} and Propositions~\ref{prop:spectral}, \ref{prop:growth}, \ref{prop:noise} and~\ref{prop:error-ratio}, which were stated without proof in Sections~\ref{sec:proof-main}--\ref{sec:verify}.

\subsection{\texorpdfstring{Proof of Lemma~\ref{lem:gaussian}}{Proof of Lemma 4.1}}\label{sec:escape-lemma}\label{subsec:proof-gaussian}

\begin{proof}[Proof of Lemma~\ref{lem:gaussian}]
The proof is partly motivated by \cite[Lemma~A.2]{CZH21}. We wish to show that, on the event $\Gamma_Y$ in \eqref{eq:Gamma-Y}, a Gaussian random variable with variance at least $\frac q{4\Lambda}$ lies in a deterministic interval of length $\frac12(q/\Lambda)^{1/2}$. Use It\^o's formula $d(e^{-\Lambda t}Y)=e^{-\Lambda t}(dY-\Lambda Y\,dt)$, and obtain:
\begin{equation}\label{eq:scalar-integrating-factor}
    e^{-\Lambda t}Y(t)=Y(0)+\sqrt q\int_0^t e^{-\Lambda s}\,dW_s
    +\int_0^t e^{-\Lambda s}\epsilon(s)\,ds .
\end{equation}

Averaging \eqref{eq:scalar-integrating-factor} over $t\in[\tau/2,\tau]$ leads to
\[
    Y(0)+Z=\frac2\tau\int_{\tau/2}^\tau e^{-\Lambda t}Y(t)\,dt
    -\frac2\tau\int_{\tau/2}^\tau\int_0^t e^{-\Lambda s}\epsilon(s)\,ds\,dt,
\]
where $Z$ is the random variable
\begin{equation}\label{eq:Z-def}
    Z\coloneqq\frac{2\sqrt q}{\tau}\int_{\tau/2}^\tau\int_0^t e^{-\Lambda s}\,dW_s\,dt .
\end{equation}
Since $e^{-\Lambda s}$ is deterministic, integration by parts, which is It\^o's formula $d(e^{-\Lambda t}W_t)=-\Lambda e^{-\Lambda t}W_t\,dt+e^{-\Lambda t}\,dW_t$, yields
\[
    \int_0^te^{-\Lambda s}\,dW_s=e^{-\Lambda t}W_t+\Lambda\int_0^te^{-\Lambda s}W_s\,ds,
\]
and consequently
\[
    Z=\frac{2\sqrt q}{\tau}\int_{\tau/2}^\tau\Bigl(e^{-\Lambda t}W_t+\Lambda\int_0^te^{-\Lambda s}W_s\,ds\Bigr)\,dt
\]
is a linear functional of the Brownian path $s\mapsto W_s$. Hence, $Z$ is a centered Gaussian random variable.

If $\tau^{-1}\int_0^\tau Y(t)^2\,dt\leq a^2$, then the Jensen inequality implies
\[
\begin{aligned}
    \Bigl(\frac2\tau\int_{\tau/2}^\tau e^{-\Lambda t}Y(t)\,dt\Bigr)^2
    &\leq\frac2\tau\int_{\tau/2}^\tau e^{-2\Lambda t}Y(t)^2\,dt\\
    &\leq\frac{2e^{-\Lambda\tau}}\tau\int_0^\tau Y(t)^2\,dt
    \leq2a^2e^{-\Lambda\tau}.
\end{aligned}
\]

Also, since $e^{-\Lambda s}\leq1$,
\[
    \left|\frac2\tau\int_{\tau/2}^\tau\int_0^t e^{-\Lambda s}\epsilon(s)\,ds\,dt\right|
    \leq\frac2\tau\int_{\tau/2}^\tau\int_0^t\abs{\epsilon(s)}\,ds\,dt
    \leq\int_0^\tau\abs{\epsilon(s)}\,ds,
\]
where in the second inequality we extended the inner integral to $[0,\tau]$ and integrated in $t$. Thus, on $\Gamma_Y$,
\begin{equation}\label{eq:Y0Z-bound}
    \abs{Y(0)+Z}\leq\sqrt2\,a\,e^{-\Lambda\tau/2}+\int_0^\tau\abs{\epsilon(s)}\,ds
    \leq\frac14\Bigl(\frac q\Lambda\Bigr)^{1/2},
\end{equation}
because $e^{\Lambda\tau}\geq128\,a^2\Lambda/q$ in \eqref{cond:window} is equivalent to $\sqrt2\,a\,e^{-\Lambda\tau/2}\leq\frac18(q/\Lambda)^{1/2}$.

For $t\in[\tau/2,\tau]$,
\begin{equation}\label{eq:split-integral}
    \int_0^t e^{-\Lambda s}\,dW_s
    =\int_0^{\tau/2}e^{-\Lambda s}\,dW_s
    +\int_{\tau/2}^t e^{-\Lambda s}\,dW_s .
\end{equation}
The integral $\int_0^{\tau/2}e^{-\Lambda s}\,dW_s$ does not depend on $t$, and therefore substituting \eqref{eq:split-integral} into the definition \eqref{eq:Z-def} of $Z$ and integrating over $t$, with $\frac2\tau\int_{\tau/2}^\tau dt=1$, yields
\[
    Z=\sqrt q\int_0^{\tau/2}e^{-\Lambda s}\,dW_s
    +\frac{2\sqrt q}{\tau}\int_{\tau/2}^\tau\int_{\tau/2}^t e^{-\Lambda s}\,dW_s\,dt .
\]
The two terms are independent because they depend on Brownian increments on $[0,\tau/2]$ and $[\tau/2,\tau]$, respectively. Their variances add, and hence, by the It\^o isometry and $\Lambda\tau\geq1$,
\[
    \mathrm{Var}\,Z\geq q\int_0^{\tau/2}e^{-2\Lambda s}\,ds
    =\frac q{2\Lambda}(1-e^{-\Lambda\tau})\geq\frac q{4\Lambda}.
\]

Since $Z$ is centered Gaussian, its probability density $f_Z:\R\to[0,\infty)$ is
\begin{equation}\label{eq:density-bound}
    f_Z(\xi)\coloneqq\frac{1}{\sqrt{2\pi\,\mathrm{Var}\,Z}}
    \exp\!\left(-\frac{\xi^2}{2\,\mathrm{Var}\,Z}\right)
    \leq\frac{1}{\sqrt{2\pi\,\mathrm{Var}\,Z}}
    \leq\Bigl(\frac\pi2\frac q\Lambda\Bigr)^{-1/2}
    \quad\text{for every }\xi\in\mathbb R.
\end{equation}
Integrating \eqref{eq:density-bound} over the interval with endpoints $-Y(0)\pm\frac14(q/\Lambda)^{1/2}$, which are deterministic, and using \eqref{eq:Y0Z-bound} for the first inequality, yields
\[
\begin{aligned}
    &\mathbb P\!\left(\frac1\tau\int_0^\tau Y(t)^2\,dt\leq a^2,\quad
    \int_0^\tau\abs{\epsilon(t)}\,dt\leq\frac18\Bigl(\frac q\Lambda\Bigr)^{1/2}\right)\\
    &\qquad\leq\mathbb P\!\left(\abs{Y(0)+Z}\leq\frac14\Bigl(\frac q\Lambda\Bigr)^{1/2}\right)
    =\int_{-Y(0)-\frac14\sqrt{q/\Lambda}}^{-Y(0)+\frac14\sqrt{q/\Lambda}}f_Z(\xi)\,d\xi\\
    &\qquad\leq\frac12\Bigl(\frac q\Lambda\Bigr)^{1/2}
    \Bigl(\frac\pi2\frac q\Lambda\Bigr)^{-1/2}
    =\frac1{\sqrt{2\pi}}\leq\frac12 .
\end{aligned}\qedhere
\]
\end{proof}

\subsection{Proof of Proposition~\ref{prop:spectral}}\label{subsec:proof-spectral}

Recall from Proposition~\ref{prop:spectral} that $\mathcal H_K$ is the completion of the space $\mathcal H$ in \eqref{eq:H-def} for the inner product $(\cdot,\cdot)_K$ in \eqref{eq:K-inner}. We begin with estimates on $Kv$ for $v\in\mathcal H$, which the proof of Proposition~\ref{prop:spectral} uses.

By \eqref{eq:fourier-operators} with $k=k_0=2\pi/L$, for $w=\cos(k_0x)g(z)\in\mathcal H$, $-\Delta w=\cos(k_0x)(k_0^2-\partial_z^2)g$, and $K=k_0^2\Delta^{-2}$ on $\mathcal H$, as in \eqref{eq:K-on-H}. For $w=\cos(k_0x)g(z)\in\mathcal H$,
\begin{equation}\label{eq:H-norms}
    \norm{w}_{L^2(D)}=(L/2)^{1/2}\norm{g}_{L^2(0,1)},\qquad\norm{w}_{L^\infty(D)}=\norm{g}_{L^\infty(0,1)}.
\end{equation}
The identities \eqref{eq:H-norms} also hold with $\partial_zw$ and $\partial_zg$, or $\partial_z^2w$ and $\partial_z^2g$, in place of $w$ and $g$, because $\partial_zw=\cos(k_0x)\partial_zg$ and $\partial_z^2w=\cos(k_0x)\partial_z^2g$.

\begin{lemma}\label{lem:vertical} For $v\in\mathcal H$,
\begin{align}
    &(v,v)_K=\frac{\norm{\Delta(Kv)}_{L^2(D)}^2}{k_0^2},\label{eq:vK-identity}\\
    &\norm{Kv}_{L^\infty(D)}+\norm{\partial_z(Kv)}_{L^\infty(D)}+\norm{\partial_z^2(Kv)}_{L^2(D)}\leq C\,(v,v)_K^{1/2}.\label{eq:Kv-bounds}
\end{align}
\end{lemma}
\begin{proof}
Observe that $\Delta(Kv)=k_0^2\Delta^{-1}v$, and consequently
\[
    \norm{\Delta(Kv)}_{L^2(D)}^2=k_0^4(\Delta^{-1}v,\Delta^{-1}v)=k_0^2(v,v)_K.
\]
Denote $Kv=\cos(k_0x)g(z)$, where $g\in H^2(0,1)$ and $g(0)=g(1)=0$. Integration by parts yields
\[
    \int_0^1g\,\partial_z^2g\,dz=-\int_0^1(\partial_zg)^2\,dz,
\]
and therefore
\[
    \norm{(k_0^2-\partial_z^2)g}_{L^2(0,1)}^2=k_0^4\norm{g}_{L^2(0,1)}^2+2k_0^2\norm{\partial_zg}_{L^2(0,1)}^2+\norm{\partial_z^2g}_{L^2(0,1)}^2.
\]
Therefore,
\begin{equation}\label{eq:d2g-bound}
    \norm{\partial_z^2g}_{L^2(0,1)}\leq\norm{(k_0^2-\partial_z^2)g}_{L^2(0,1)},
\end{equation}
and, by the Sobolev embedding $H^2(0,1)\subset C^1[0,1]$,
\begin{equation}\label{eq:g-sup-bound}
    \norm{g}_{L^\infty(0,1)}+\norm{\partial_zg}_{L^\infty(0,1)}\leq C\norm{(k_0^2-\partial_z^2)g}_{L^2(0,1)}.
\end{equation}
The estimate \eqref{eq:Kv-bounds} follows from \eqref{eq:d2g-bound}, \eqref{eq:g-sup-bound}, \eqref{eq:H-norms} with $w=Kv$ and with $w=\Delta(Kv)=-\cos(k_0x)(k_0^2-\partial_z^2)g$, and \eqref{eq:vK-identity}.
\end{proof}

For $\rho$ as in Proposition~\ref{prop:spectral} and $v,w\in\mathcal H$, the Cauchy--Schwarz inequality in $x$ and then in $z$, and Lemma~\ref{lem:vertical} with $\norm{Kw}_{L^2(D)}\leq L^{1/2}\norm{Kw}_{L^\infty(D)}$ imply
\begin{equation}\label{eq:bounded-part}
\begin{aligned}
    \Bigl|\Ray\bigl((1-\partial_z\rho)Kv,Kw\bigr)\Bigr|
    &\leq\Ray\int_0^1\abs{1-\partial_z\rho}\Bigl(\int_0^L\abs{Kv\,Kw}\,dx\Bigr)dz\\
    &\leq L^{1/2}\Ray\,\norm{1-\partial_z\rho}_{L^2(0,1)}\norm{Kv}_{L^\infty(D)}\norm{Kw}_{L^2(D)}\\
    &\leq C\Ray\,\norm{1-\partial_z\rho}_{L^2(0,1)}(v,v)_K^{1/2}(w,w)_K^{1/2}.
\end{aligned}
\end{equation}

\begin{proof}[Proof of Proposition~\ref{prop:spectral}]
By the definition of $\mathcal A_\rho$ in \eqref{eq:A-rho},
\begin{equation}\label{eq:A-split}
    \mathcal A_\rho=\Delta+\Ray(1-\partial_z\rho)K.
\end{equation}
We show that the summand $\Ray(1-\partial_z\rho)K$ is a bounded symmetric operator on $\mathcal H_K$, and that the summand $\Delta$ is self-adjoint on $\mathcal H_K$, with domain $\mathcal H$ and compact resolvent. The Kato--Rellich theorem then yields the desired properties of the sum $\mathcal A_\rho$.

First, we consider $\Ray(1-\partial_z\rho)K$. For $v,w\in\mathcal H$, since $(v,w)_K=(Kv,w)$ and $K$ is symmetric for $(\cdot,\cdot)$,
\[
    \bigl(\Ray(1-\partial_z\rho)Kv,w\bigr)_K=\Ray\bigl((1-\partial_z\rho)Kv,Kw\bigr),
\]
which is symmetric in $v$ and $w$. By \eqref{eq:bounded-part},
\[
    \Bigl|\bigl(\Ray(1-\partial_z\rho)Kv,w\bigr)_K\Bigr|\leq C\Ray\,\norm{1-\partial_z\rho}_{L^2(0,1)}(v,v)_K^{1/2}(w,w)_K^{1/2}.
\]
Therefore, $\Ray(1-\partial_z\rho)K$ extends to a bounded symmetric operator on $\mathcal H_K$, by density.

Next, we consider $\Delta$. By \eqref{eq:fourier-operators}, the functions $c_{k_0,n}$ are orthogonal for $(\cdot,\cdot)_K$, with $(c_{k_0,n},c_{k_0,n})_K=\frac{k_0^2}{(k_0^2+n^2\pi^2)^2}$, and their span is dense in $\mathcal H_K$. Therefore, every $v\in\mathcal H_K$ is a series that converges in $\mathcal H_K$,
\[
    v=\sum_na_nc_{k_0,n},
\]
and, for those $v$ for which the following series converges in $\mathcal H_K$, we define
\[
    \Delta v\coloneqq-\sum_n(k_0^2+n^2\pi^2)\,a_nc_{k_0,n}.
\]
For $v\in\mathcal H\cap H^1_0(D)$ with $\Delta v\in L^2(D)$, this series is the Laplacian of $v$, by \eqref{eq:fourier-operators}. The operator $\Delta$ is diagonal in the $(\cdot,\cdot)_K$-orthogonal basis $\{c_{k_0,n}\}$, with eigenvalues $-(k_0^2+n^2\pi^2)\to-\infty$, and hence $\Delta$ is self-adjoint with compact resolvent. The domain of $\Delta$ is $\mathcal H$, because
\[
    (\Delta v,\Delta v)_K=\sum_n(k_0^2+n^2\pi^2)^2a_n^2\,(c_{k_0,n},c_{k_0,n})_K=k_0^2\sum_na_n^2,
\]
which is finite exactly when $v\in\mathcal H$. Moreover, since every eigenvalue of $\Delta$ is at most $-(k_0^2+\pi^2)$,
\begin{equation}\label{eq:Delta-upper}
    (\Delta v,v)_K\leq-(k_0^2+\pi^2)(v,v)_K,\qquad v\in\mathcal H.
\end{equation}

Finally, the Kato--Rellich theorem \cite[Theorem~X.12]{RS2} states that a self-adjoint operator plus a bounded symmetric operator is self-adjoint on the same domain. Thus, by \eqref{eq:A-split}, $\mathcal A_\rho$ is self-adjoint on $\mathcal H_K$ with domain $\mathcal H$. For $\mu\in\mathbb C\setminus\R$, the resolvent identity
\[
    (\mathcal A_\rho-\mu)^{-1}=(\Delta-\mu)^{-1}-(\Delta-\mu)^{-1}\,\Ray(1-\partial_z\rho)K\,(\mathcal A_\rho-\mu)^{-1}
\]
shows that $\mathcal A_\rho$ has compact resolvent, because $(\Delta-\mu)^{-1}$ is compact. Also, $\mathcal A_\rho$ is bounded above, by \eqref{eq:Delta-upper} and because $\Ray(1-\partial_z\rho)K$ is bounded.
\end{proof}

\subsection{Proof of Proposition~\ref{prop:growth}}\label{subsec:proof-growth}

\begin{proof} Since the largest eigenvalue of $\mathcal A_\rho$ is given by \eqref{eq:largest-eigenvalue}, it suffices to find one $v=c_{k_0,n}$ with $(\mathcal A_\rho v,v)_K\geq\lambda_*(v,v)_K$. We first claim that there is $n\leq1+2\beta$ with
\begin{equation}\label{eq:n-choice}
    1+\int_0^1\partial_z\rho\cos(2n\pi z)\,dz\geq\delta,
\end{equation}
where $\delta\coloneqq1-\frac1{\sqrt2}$. Indeed, the functions $\sqrt2\cos(2n\pi z)$ are orthonormal, and therefore Bessel's inequality, with the factor $\frac12$ from the normalization $\sqrt2$, and the assumption $\int_0^1(\partial_z\rho)^2dz<2\beta$ imply
\[
    \sum_n\Bigl(\int_0^1\partial_z\rho\cos(2n\pi z)\,dz\Bigr)^2\leq\frac12\int_0^1(\partial_z\rho)^2dz<\beta .
\]
If $\int_0^1\partial_z\rho\cos(2n\pi z)\,dz<-\frac1{\sqrt2}$ for every $n\leq\lceil2\beta\rceil$, the left side would exceed $\frac12\lceil2\beta\rceil\geq\beta$, a contradiction.

Denote $v\coloneqq c_{k_0,n}$ for $n$ satisfying $1\leq n\leq1+2\beta$ and \eqref{eq:n-choice}. By \eqref{eq:fourier-operators}, $Kv=(v,v)_K\,v$, with $(v,v)_K=\frac{k_0^2}{(k_0^2+n^2\pi^2)^2}$. Since
\[
    \int_0^Lc_{k_0,n}^2\,dx=2\sin^2(n\pi z)=1-\cos(2n\pi z),
\]
\eqref{eq:A-form} yields
\[
    (\mathcal A_\rho v,v)_K=-(k_0^2+n^2\pi^2)(v,v)_K+\Ray\,(v,v)_K^2\int_0^1(1-\partial_z\rho)\bigl(1-\cos(2n\pi z)\bigr)\,dz.
\]
Note that
\[
    \int_0^1\cos(2n\pi z)\,dz=0,\qquad\int_0^1\partial_z\rho\,dz=0,
\]
the second because $\rho\in H^1_0(0,1)$. Dividing by $(v,v)_K$, using \eqref{eq:n-choice}, we obtain
\[
\begin{aligned}
\frac{(\mathcal A_\rho v,v)_K}{(v,v)_K}
&=\frac{\Ray\,k_0^2}{(k_0^2+n^2\pi^2)^2}
  \Bigl(1+\int_0^1\partial_z\rho\cos(2n\pi z)\,dz\Bigr)
  -(k_0^2+n^2\pi^2)\\
&\geq\frac{\delta\,\Ray\,k_0^2}{(k_0^2+n^2\pi^2)^2}-(k_0^2+n^2\pi^2).
\end{aligned}
\]
Since $n\leq1+2\beta$, $k_0^2+n^2\pi^2\leq(k_0^2+\pi^2)(1+2\beta)^2$, and consequently the assumption $\Ray\geq C(1+2\beta)^6$, with $C\geq\frac2\delta(k_0^2+\pi^2)^3/k_0^2$, which depends only on $k_0=2\pi/L$, implies
\[
    k_0^2+n^2\pi^2\leq\frac{\delta\,\Ray\,k_0^2}{2(k_0^2+n^2\pi^2)^2}.
\]
Therefore,
\[
\begin{aligned}
\sup_{w\in\mathcal H,\,w\neq0}\frac{(\mathcal A_\rho w,w)_K}{(w,w)_K}
&\geq\frac{(\mathcal A_\rho v,v)_K}{(v,v)_K}\\
&\geq\frac{\delta\,\Ray\,k_0^2}{2(k_0^2+n^2\pi^2)^2}\geq\frac{\delta\,\Ray\,k_0^2}{2(1+2\beta)^4(k_0^2+\pi^2)^2}\\
&\geq\lambda_*=c\,\Ray\,(1+2\beta)^{-4}\quad\text{for }c\leq\frac{\delta\,k_0^2}{2(k_0^2+\pi^2)^2}. \end{aligned}\qedhere
\]
\end{proof}

\subsection{Proof of Proposition~\ref{prop:noise}}\label{subsec:proof-noise}

\begin{proof}
We wish to show that $q/\Lambda$ has a lower bound that is independent of $\theta_0\in\mathcal S_{Q<\beta}$. Recall the formula \eqref{eq:q-formula} for $q$. The Cauchy--Schwarz inequality for $(\cdot,\cdot)_K$ and $(\Psi,\Psi)_K=1$ imply
\[
    \norm{K\Psi}_{L^2(D)}^2=(K\Psi,\Psi)_K\leq(K\Psi,K\Psi)_K^{1/2},
\]
and therefore
\begin{equation}\label{eq:noise-profile-bound}
q\geq\frac L4C_L^2\sigma^2\norm {{K\Psi}}_{L^2(D)}^4 .
\end{equation}

Next, we estimate $\norm{K\Psi}_{L^\infty(D)}$. By \eqref{eq:Kv-bounds} and $(\Psi,\Psi)_K=1$, $\norm{\partial_z^2(K\Psi)}_{L^2(D)}\leq C$. Integration by parts in $z$ yields
\begin{equation}\label{eq:dzKPsi}
    \norm{\partial_z(K\Psi)}_{L^2(D)}^2=-\bigl(K\Psi,\partial_z^2(K\Psi)\bigr)\leq C\norm{K\Psi}_{L^2(D)}.
\end{equation}
For every $g\in H^1_0(0,1)$,
\begin{equation}\label{eq:g-sup}
    g(z)^2=2\int_0^zg\,\partial_zg\,dz'\leq2\norm g_{L^2(0,1)}\norm{\partial_zg}_{L^2(0,1)}\qquad\text{for any }z\in(0,1).
\end{equation}
The identities \eqref{eq:H-norms}, the estimate \eqref{eq:g-sup} with $K\Psi=\cos(k_0x)g(z)$, and then \eqref{eq:dzKPsi} imply
\begin{equation}\label{eq:KPsi-sup}
    \norm{K\Psi}_{L^\infty(D)}\leq C\norm{K\Psi}_{L^2(D)}^{1/2}\norm{\partial_z(K\Psi)}_{L^2(D)}^{1/2}\leq C\norm{K\Psi}_{L^2(D)}^{3/4}.
\end{equation}

Since $\mathcal A_{\bar\theta_0}\Psi=\Lambda\Psi$ and $(\Psi,\Psi)_K=1$, $\Lambda=(\mathcal A_{\bar\theta_0}\Psi,\Psi)_K$, and the first term $k_0^2(\Delta^{-1}\Psi,\Psi)$ of $(\mathcal A_{\bar\theta_0}\Psi,\Psi)_K$ in \eqref{eq:A-form} is at most $0$, by \eqref{eq:fourier-operators}. Also, since $\bar\theta_0\in H^1_0(0,1)$, and by \eqref{eq:low-set-bounds},
\[
    \norm{1-\partial_z\bar\theta_0}_{L^2(0,1)}^2=1-2\int_0^1\partial_z\bar\theta_0\,dz+\int_0^1(\partial_z\bar\theta_0)^2dz=1+\int_0^1(\partial_z\bar\theta_0)^2dz<1+2\beta.
\]
By the Cauchy--Schwarz inequality and \eqref{eq:KPsi-sup},
\begin{equation}\label{eq:growth-profile-bound}
\begin{aligned}
    \Lambda&=(\mathcal A_{\bar\theta_0}\Psi,\Psi)_K
    \leq\Ray\bigl((1-\partial_z\bar\theta_0)K\Psi,K\Psi\bigr)\\
    &\leq C\Ray\,\norm{1-\partial_z\bar\theta_{0}}_{L^2(0,1)}\,\norm {{K\Psi}}_{L^\infty(D)}\norm {{K\Psi}}_{L^2(D)}
    \leq C\Ray\,(1+2\beta)^{1/2}\norm {{K\Psi}}_{L^2(D)}^{7/4}.
\end{aligned}
\end{equation}

Then \eqref{eq:construct-assumptions}, \eqref{eq:low-set-bounds} and Proposition~\ref{prop:growth} imply
\[
    c\,\Ray\,(1+2\beta)^{-4}=\lambda_*\leq\Lambda\leq C\Ray\,(1+2\beta)^{1/2}\norm{K\Psi}_{L^2(D)}^{7/4},
\]
and consequently $\norm{K\Psi}_{L^2(D)}\geq c\,(1+2\beta)^{-18/7}$. Dividing \eqref{eq:noise-profile-bound} by \eqref{eq:growth-profile-bound} and using this lower bound on $\norm{K\Psi}_{L^2(D)}$ yields
\[
\frac q\Lambda\geq\frac{c\,\sigma^2}\Ray\,(1+2\beta)^{-1/2}\norm {{K\Psi}}_{L^2(D)}^{9/4}\geq\frac{c\,\sigma^2}\Ray\,(1+2\beta)^{-44/7}. \qedhere
\]
\end{proof}

\subsection{Proof of Proposition~\ref{prop:error-ratio}}\label{subsec:proof-errors}

\begin{proof}
Let $\beta$, $\Ray$, $\theta_0$, $\theta$, $\Lambda$, $\Psi$, $q$, $\epsilon_1$, $\epsilon_2$ and $\tau$ be as in Proposition~\ref{prop:error-ratio}. We wish to show \eqref{eq:eps-integrable} almost surely, and \eqref{eq:error-ratio} almost surely on $\Gamma_Q$. Recall from \eqref{eq:mean-def} that $\bar\theta$ is the horizontal mean of $\theta$ and $\theta'=\theta-\bar\theta$. We first collect estimates at a fixed time $t$ with $\theta(t)\in H^1_0(D)$, where $Q=Q(\theta(t))$. The two terms $\frac\Ray{2L}(K\theta',\theta')$ and $\frac1{2L}\norm{\nabla\theta'}_{L^2(D)}^2$ of $Q$ in \eqref{eq:mean-split} yield
\begin{equation}\label{eq:theta-parts}
    (K\theta',\theta')\leq\frac{2L}{\Ray}Q,
    \qquad \norm{\nabla\theta'}_{L^2(D)}^2\leq2LQ.
\end{equation}
Also, by \eqref{eq:stokes-energy} and $Q\geq\frac1{2L}(u_2,\theta)$ from \eqref{eq:Q-def},
\begin{equation}\label{eq:grad-u}
    \norm{\nabla\vu}_{L^2(D)}^2=\Ray\,(u_2,\theta)\leq2L\Ray Q.
\end{equation}
Next, by interpolating between the two estimates in \eqref{eq:theta-parts}, we obtain estimates on $\theta'$ in $L^2(D)$ as follows. Since $(\theta',e_n)=0$ for every $n$, Parseval's identity and \eqref{eq:H1-fourier} yield
\[
\begin{aligned}
    &\norm{\theta'}_{L^2(D)}^2=\sum_{k,n}\bigl[(\theta',c_{k,n})^2+(\theta',s_{k,n})^2\bigr],\\
    &\norm{\nabla\theta'}_{L^2(D)}^2=\sum_{k,n}(k^2+n^2\pi^2)\bigl[(\theta',c_{k,n})^2+(\theta',s_{k,n})^2\bigr],
\end{aligned}
\]
and, since $k\geq k_0$ in each sum, \eqref{eq:fourier-operators} implies
\[
    (K\theta',\theta')=\sum_{k,n}\frac{k^2}{(k^2+n^2\pi^2)^2}\bigl[(\theta',c_{k,n})^2+(\theta',s_{k,n})^2\bigr]\geq k_0^2\sum_{k,n}\frac{\bigl[(\theta',c_{k,n})^2+(\theta',s_{k,n})^2\bigr]}{(k^2+n^2\pi^2)^2}.
\]
Then H\"older's inequality with exponents $3$ and $\frac32$ and \eqref{eq:theta-parts} imply
\begin{equation}\label{eq:theta-bound}
\begin{aligned}
    \norm{\theta'}_{L^2(D)}^2&=\sum_{k,n}\Bigl(\frac{\bigl[(\theta',c_{k,n})^2+(\theta',s_{k,n})^2\bigr]}{(k^2+n^2\pi^2)^2}\Bigr)^{1/3}\Bigl((k^2+n^2\pi^2)\bigl[(\theta',c_{k,n})^2+(\theta',s_{k,n})^2\bigr]\Bigr)^{2/3}\\
    &\leq\Bigl(\sum_{k,n}\frac{\bigl[(\theta',c_{k,n})^2+(\theta',s_{k,n})^2\bigr]}{(k^2+n^2\pi^2)^2}\Bigr)^{1/3}\Bigl(\sum_{k,n}(k^2+n^2\pi^2)\bigl[(\theta',c_{k,n})^2+(\theta',s_{k,n})^2\bigr]\Bigr)^{2/3}\\
    &\leq\Bigl(\frac{(K\theta',\theta')}{k_0^2}\Bigr)^{1/3}\norm{\nabla\theta'}_{L^2(D)}^{4/3}\\
    &\leq\Bigl(\frac{2LQ}{k_0^2\Ray}\Bigr)^{1/3}(2LQ)^{2/3}=\frac{2LQ}{(k_0^2\Ray)^{1/3}}.
\end{aligned}
\end{equation}

We begin with $\epsilon_2$. Fix $t\geq0$ with $\theta(t)\in H^1_0(D)$, and recall from \eqref{eq:eps-def} that $\epsilon_2=(\theta',f)$, where $f\coloneqq\vu(t)\cdot\nabla K\Psi$ as a function $D\to\R$. We estimate $\epsilon_2$ by the $H^{-1}(D)$--$H^1_0(D)$ pairing:
\begin{equation}\label{eq:epsilon-two-pairing}
    \abs{\epsilon_2}\leq\norm{\theta'}_{H^{-1}(D)}\norm{f}_{H^1_0(D)}.
\end{equation}
Since $(\theta',e_n)=0$ for every $n$, only the first sum in \eqref{eq:Hminus1-fourier} contributes to $\norm{\theta'}_{H^{-1}(D)}^2$. Therefore, the Cauchy--Schwarz inequality, $k\geq k_0$, \eqref{eq:theta-parts} and \eqref{eq:theta-bound} imply
\begin{equation}\label{eq:theta-Hminus1}
\begin{aligned}
    \norm{\theta'}_{H^{-1}(D)}^2&=\sum_{k,n}\frac{\bigl[(\theta',c_{k,n})^2+(\theta',s_{k,n})^2\bigr]}{k^2+n^2\pi^2}\\
    &\leq\Bigl(\sum_{k,n}\frac{\bigl[(\theta',c_{k,n})^2+(\theta',s_{k,n})^2\bigr]}{(k^2+n^2\pi^2)^2}\Bigr)^{1/2}\Bigl(\sum_{k,n}\bigl[(\theta',c_{k,n})^2+(\theta',s_{k,n})^2\bigr]\Bigr)^{1/2}\\
    &\leq\Bigl(\frac{(K\theta',\theta')}{k_0^2}\Bigr)^{1/2}\norm{\theta'}_{L^2(D)}\leq C\,\Ray^{-2/3}Q.
\end{aligned}
\end{equation}

Since $f=0$ at $z=0,1$, integration by parts yields $(\nabla f,\nabla c_{k,n})=-(f,\Delta c_{k,n})=(k^2+n^2\pi^2)(f,c_{k,n})$, and likewise for $s_{k,n}$ and $e_n$, and therefore $\norm{f}_{H^1_0(D)}=\norm{\nabla f}_{L^2(D)}$. The product rule yields
\[
    \nabla f=\nabla u_1\,\partial_x(K\Psi)+\nabla u_2\,\partial_z(K\Psi)+u_1\nabla\partial_x(K\Psi)+u_2\nabla\partial_z(K\Psi).
\]
Since $K\Psi\in\mathcal H$ (see \eqref{eq:H-def}), then $K\Psi=\cos(k_0x)g(z)$ for some $g\in H^2(0,1)$. Lemma~\ref{lem:vertical}, \eqref{eq:H-norms} and $(\Psi,\Psi)_K=1$ imply
\[
\begin{gathered}
    \norm{\nabla K\Psi}_{L^\infty(D)}\leq k_0\norm{K\Psi}_{L^\infty(D)}+\norm{\partial_z(K\Psi)}_{L^\infty(D)}\leq C,\\
    \norm{\partial_x^2(K\Psi)}_{L^\infty(D)}+\norm{\partial_x\partial_z(K\Psi)}_{L^\infty(D)}\leq k_0^2\norm{K\Psi}_{L^\infty(D)}+k_0\norm{\partial_z(K\Psi)}_{L^\infty(D)}\leq C,\\
    \norm{\partial_z^2g}_{L^2(0,1)}\leq C.
\end{gathered}
\]
Since $u_2=0$ at $z=0,1$ and $u_1$ has mean $0$, Poincar\'e's inequality implies $\norm\vu_{L^2(D)}\leq C\norm{\nabla\vu}_{L^2(D)}$.
Hence, the terms of $\nabla f$ except $u_2\partial_z^2(K\Psi)$ can be estimated as
\[
\begin{aligned}
    \norm{\nabla u_1\,\partial_x(K\Psi)+\nabla u_2\,\partial_z(K\Psi)}_{L^2(D)}&\leq2\norm{\nabla\vu}_{L^2(D)}\norm{\nabla K\Psi}_{L^\infty(D)}\leq C\norm{\nabla\vu}_{L^2(D)},\\
    \norm{u_1\nabla\partial_x(K\Psi)}_{L^2(D)}+\norm{u_2\partial_x\partial_z(K\Psi)}_{L^2(D)}&\leq C\norm\vu_{L^2(D)}\leq C\norm{\nabla\vu}_{L^2(D)}.
\end{aligned}
\]
To treat $u_2\partial_z^2(K\Psi)$ we note that, since $u_2=0$ at $z=0$, the Cauchy--Schwarz inequality in $\xi$ implies, for every $z\in[0,1]$,
\begin{equation}
\begin{aligned}
    \int_0^Lu_2(x,z)^2\,dx&=\int_0^L\Bigl(\int_0^z\partial_zu_2(x,\xi)\,d\xi\Bigr)^2dx\leq\int_0^L\int_0^z\bigl(\partial_zu_2(x,\xi)\bigr)^2\,d\xi\,dx\\
    &\leq\norm{\partial_zu_2}_{L^2(D)}^2 .
\end{aligned}\label{eq:u2-height}
\end{equation}
Thus, from $\abs{\cos(k_0x)}\leq1$ follows
\[
    \norm{u_2\partial_z^2(K\Psi)}_{L^2(D)}^2\leq\int_0^1(\partial_z^2g)^2\int_0^Lu_2^2\,dx\,dz\leq C\norm{\partial_zu_2}_{L^2(D)}^2.
\]
Altogether, with \eqref{eq:grad-u},
\begin{equation}\label{eq:grad-f-bound}
    \norm{\nabla f}_{L^2(D)}\leq C\norm{\nabla\vu}_{L^2(D)}\leq C(\Ray Q)^{1/2}.
\end{equation}
Substituting \eqref{eq:theta-Hminus1} and \eqref{eq:grad-f-bound} into \eqref{eq:epsilon-two-pairing}, we obtain
\begin{equation}\label{eq:error-two-pointwise}
    \abs{\epsilon_2(t)}\leq C\,\Ray^{-1/3}Q^{1/2}\cdot\Ray^{1/2}Q^{1/2}=C\,\Ray^{1/6}Q(\theta(t)),
\end{equation}
and hence on the event $\Gamma_Q$ in \eqref{eq:Gamma-Q},
\begin{equation}\label{eq:error-two}
    \int_0^\tau\abs{\epsilon_2}\,dt\leq C\Ray^{1/6}\beta\tau .
\end{equation}

Next we turn our attention to $\epsilon_1$. To estimate the difference $\partial_z\bar\theta_0-\partial_z\bar\theta(t)$ in $\epsilon_1$, we derive an equation for $\bar\theta$. First, define $J:[0,\infty)\times[0,1]\to\R$ as the horizontal average of $u_2\theta'$, which appears in this equation:
\[
    J(r,z)\coloneqq\frac1L\int_0^L u_2(r,x,z)\theta'(r,x,z)\,dx .
\]
Let $\varphi\in H^1_0(0,1)$. Since the noise coefficients $\zeta_n$ in \eqref{eq:noise}, for each $n$, and $u_2$ have zero average in $x$, $(\zeta_n,\varphi)=0$ and $\Ray(K\theta,\varphi)=(u_2,\varphi)=0$, and therefore \eqref{eq:weak-solution} with $\vartheta=\varphi$ becomes
\[
    (\theta(t),\varphi)=(\theta_0,\varphi)+\int_0^t\Bigl(-(\nabla\theta(r),\nabla\varphi)+(\theta(r),\vu(r)\cdot\nabla\varphi)\Bigr)\,dr\qquad\text{for all }t\geq0,
\]
where, since $\varphi$ is independent of $x$, by Fubini's theorem,
\[
\begin{aligned}
    (\theta,\varphi)&=L\int_0^1\bar\theta\varphi\,dz,\\
    (\nabla\theta,\nabla\varphi)&=L\int_0^1\partial_z\bar\theta\,\partial_z\varphi\,dz,\\
    (\theta,\vu\cdot\nabla\varphi)=(\theta,u_2\partial_z\varphi)&=L\int_0^1J\,\partial_z\varphi\,dz.
\end{aligned}
\]
Dividing by $L$, we obtain
\begin{equation}\label{eq:mean-weak}
    \int_0^1\bigl(\bar\theta(t)-\bar\theta_0\bigr)\varphi\,dz=-\int_0^1\Bigl(\int_0^t\bigl(\partial_z\bar\theta(r)-J(r)\bigr)\,dr\Bigr)\partial_z\varphi\,dz,
\end{equation}
almost surely, for all $t\geq0$ and every $\varphi\in H^1_0(0,1)$, that is, $\partial_t\bar\theta=\partial_z(\partial_z\bar\theta-J)$ weakly.

Fix $t$ with $\theta(t)\in H^1_0(D)$. Since $\partial_z\bar\theta_0-\partial_z\bar\theta(t)$ depends only on $z$, the definition \eqref{eq:eps-def} of $\epsilon_1$ yields
\[
    \epsilon_1(t)=\Ray\int_0^1(\partial_z\bar\theta_0-\partial_z\bar\theta(t))\eta_t\,dz,
\]
where we denote $\eta_t:[0,1]\to\R$ by
\[
    \eta_t(z)\coloneqq\int_0^LK(P_{\mathcal H}\theta(t))\,K\Psi\,dx.
\]
Then $\eta_t=0$ and $\partial_z\eta_t=0$ at $z=0,1$, since $KP_{\mathcal H}\theta(t)=0$ and $K\Psi=0$ for $z=0,1$. Integrating by parts in $\epsilon_1$ and using \eqref{eq:mean-weak} with $\varphi=\partial_z\eta_t$, we obtain
\begin{equation}\label{eq:error-one-representation}
\begin{aligned}
\epsilon_1(t)&=\Ray\int_0^1\bigl(\bar\theta(t,z)-\bar\theta_0(z)\bigr)\partial_z\eta_t(z)\,dz\\
&=-\Ray\int_0^1\Bigl(\int_0^t\bigl(\partial_z\bar\theta(r,z)-J(r,z)\bigr)\,dr\Bigr)
  \partial_z^2\eta_t(z)\,dz .
\end{aligned}
\end{equation}
Differentiating $\eta_t$ twice under the integral, we obtain
\[
\begin{aligned}
    \partial_z^2\eta_t=\int_0^L\Bigl(&\partial_z^2(KP_{\mathcal H}\theta)\,K\Psi+2\,\partial_z(KP_{\mathcal H}\theta)\,\partial_z(K\Psi)\\
    &+KP_{\mathcal H}\theta\,\partial_z^2(K\Psi)\Bigr)\,dx,
\end{aligned}
\]
and the Cauchy--Schwarz inequality in $x$ implies
\[
\begin{aligned}
    \norm{\partial_z^2\eta_t}_{L^2(0,1)}\leq L^{1/2}\Bigl(&\norm{\partial_z^2(KP_{\mathcal H}\theta)}_{L^2(D)}\norm{K\Psi}_{L^\infty(D)}+2\norm{\partial_z(KP_{\mathcal H}\theta)}_{L^\infty(D)}\norm{\partial_z(K\Psi)}_{L^2(D)}\\
    &+\norm{KP_{\mathcal H}\theta}_{L^\infty(D)}\norm{\partial_z^2(K\Psi)}_{L^2(D)}\Bigr).
\end{aligned}
\]
Lemma~\ref{lem:vertical}, $(\Psi,\Psi)_K=1$ and the estimate $(P_{\mathcal H}\theta,P_{\mathcal H}\theta)_K\leq\frac{2L}\Ray Q(\theta(t))$ from \eqref{eq:coordinate-bound} imply
\begin{equation}\label{eq:g-bound}
    \norm{\partial_z^2\eta_t}_{L^2(0,1)}\leq C(P_{\mathcal H}\theta,P_{\mathcal H}\theta)_K^{1/2}\leq C\Bigl(\frac{Q(\theta(t))}\Ray\Bigr)^{1/2}.
\end{equation}

The Cauchy--Schwarz inequality in $x$, and then \eqref{eq:u2-height}, yield
\[
\begin{aligned}
\abs{J(r,z)}&=\frac1L\Bigl|\int_0^Lu_2(r,x,z)\,\theta'(r,x,z)\,dx\Bigr|\\
&\leq\frac1L\Bigl(\int_0^Lu_2(r,x,z)^2\,dx\Bigr)^{1/2}\Bigl(\int_0^L\theta'(r,x,z)^2\,dx\Bigr)^{1/2}\\
&\leq\frac1L\norm{\partial_zu_2(r)}_{L^2(D)}\Bigl(\int_0^L\theta'(r,x,z)^2\,dx\Bigr)^{1/2}.
\end{aligned}
\]
Squaring, integrating over $z\in[0,1]$, and then using \eqref{eq:grad-u} and \eqref{eq:theta-bound} at time $r$, we obtain
\begin{equation}\label{eq:J-bound}
\begin{aligned}
    \norm{J(r,\cdot)}_{L^2(0,1)}&\leq\frac1L\norm{\nabla\vu(r)}_{L^2(D)}\,\norm{\theta'(r)}_{L^2(D)}\\
    &\leq C\Ray^{1/2}Q(\theta(r))^{1/2}\,\Ray^{-1/6}Q(\theta(r))^{1/2}=C\Ray^{1/3}Q(\theta(r)).
\end{aligned}
\end{equation}
From the Cauchy--Schwarz inequality in $z$ and Minkowski's inequality, applied to \eqref{eq:error-one-representation}, follows
\[
    \abs{\epsilon_1(t)}\leq\Ray\norm{\partial_z^2\eta_t}_{L^2(0,1)}\int_0^t\Bigl(\norm{\partial_z\bar\theta(r)}_{L^2(0,1)}+\norm{J(r)}_{L^2(0,1)}\Bigr)\,dr.
\]
Using \eqref{eq:g-bound}, $\norm{\partial_z\bar\theta(r)}_{L^2(0,1)}\leq(2Q(\theta(r)))^{1/2}$ from \eqref{eq:mean-split}, and \eqref{eq:J-bound} then yields
\begin{equation}\label{eq:error-one-pointwise}
    \abs{\epsilon_1(t)}\leq C\Ray^{1/2}Q(\theta(t))^{1/2}\int_0^t\Bigl(Q(\theta(r))^{1/2}+\Ray^{1/3}Q(\theta(r))\Bigr)\,dr .
\end{equation}
Integrating over $t\in[0,\tau]$, and extending the inner integral to $[0,\tau]$, we obtain
\begin{equation}\label{eq:error-one-integrated}
    \int_0^\tau\abs{\epsilon_1}\,dt\leq C\Ray^{1/2}\Bigl(\int_0^\tau Q(\theta(t))^{1/2}\,dt\Bigr)\Bigl(\int_0^\tau Q(\theta(r))^{1/2}\,dr+\Ray^{1/3}\int_0^\tau Q(\theta(r))\,dr\Bigr).
\end{equation}
By Proposition~\ref{prop:wellposed}(i) and \eqref{eq:u2-theta-bound},
\[
    \E_{\theta_0}\int_0^\tau Q(\theta(t))\,dt<\infty,
\]
and therefore the function $t\mapsto Q(\theta(t))$ is integrable on $[0,\tau]$, almost surely. Consequently, \eqref{eq:error-two-pointwise}, \eqref{eq:error-one-integrated} and the Cauchy--Schwarz inequality in time imply \eqref{eq:eps-integrable}.

On $\Gamma_Q$, by \eqref{eq:Gamma-Q}, the Cauchy--Schwarz inequality in time implies
\[
    \int_0^\tau Q(\theta(t))^{1/2}\,dt\leq\tau^{1/2}\Bigl(\int_0^\tau Q(\theta(t))\,dt\Bigr)^{1/2}\leq(2\beta)^{1/2}\tau.
\]
Consequently, since $\beta^{1/2}\leq\Ray^{1/3}\beta$ for $\beta\geq1$ and $\Ray\geq1$, which hold by \eqref{eq:construct-assumptions}, we have
\begin{equation}\label{eq:error-one}
    \int_0^\tau\abs{\epsilon_1}\,dt\leq C\Ray^{1/2}\beta^{1/2}\tau\cdot\bigl(\beta^{1/2}+\Ray^{1/3}\beta\bigr)\tau\leq C\Ray^{5/6}\beta^{3/2}\tau^2 .
\end{equation}

Finally, by \eqref{eq:window-length} and $1+2\beta\leq3\beta$, $\tau\leq C\beta^5/\Ray$. Therefore, \eqref{eq:error-two} and \eqref{eq:error-one} imply
\[
    \int_0^\tau\abs{\epsilon_1}\,dt\leq C\beta^{23/2}\Ray^{-7/6},\qquad\int_0^\tau\abs{\epsilon_2}\,dt\leq C\beta^6\Ray^{-5/6},
\]
and Proposition~\ref{prop:noise} implies $(q/\Lambda)^{1/2}\geq c\,\sigma\,\Ray^{-1/2}\beta^{-22/7}$, and therefore \eqref{eq:error-ratio} holds.
\end{proof}

\appendix
\section{Well-posedness and the Markov property}\label{sec:wellposed}\label{subsec:proof-wellposed}

To prove Proposition~\ref{prop:wellposed}, we use the variational approach to stochastic partial differential equations (see \cite[Chapters~4 and~5]{LR15} for reference). In Section~\ref{subsec:noise-cylindrical}, we write the noise in \eqref{eq:theta-equation} as a stochastic integral with respect to a cylindrical Wiener process. In Section~\ref{subsec:variational}, we recall an existence theorem for variational solutions and It\^o's formula for the squared norm. In Section~\ref{subsec:proof-prop-wellposed}, we prove Proposition~\ref{prop:wellposed}.

\subsection{The noise in terms of a cylindrical Wiener process}\label{subsec:noise-cylindrical}

In the variational approach, the noise has the form $B\,d\mathcal W$, where $\mathcal W$ is a cylindrical Wiener process in a separable Hilbert space (see \cite[Section~2.5.1]{LR15} for reference), and $B$ takes values in the Hilbert--Schmidt operators from this space to $L^2(D)$. In this section, we write the noise $\sigma\sum_n\zeta_n\,dW^n_t$ in \eqref{eq:theta-equation} in this form.

Let $g_n\coloneqq(0,\dots,0,1,0,\dots)\in\ell^2$, with $1$ in the $n$th entry. Because the $W^n$ are independent standard Brownian motions, the formal series
\begin{equation}\label{eq:cyl-W}
    \mathcal W(t)\coloneqq\sum_nW^n_t\,g_n
\end{equation}
defines a cylindrical Wiener process $\mathcal W$ in $\ell^2$ (see \cite[Proposition~2.5.2]{LR15}). Let $B$ be the linear operator
\begin{equation}\label{eq:B-def}
    B:\ell^2\to L^2(D),\qquad Bg_n\coloneqq\sigma\zeta_n.
\end{equation}
Because $\sum_n\norm{\zeta_n}_{L^2(D)}^2=1$ by \eqref{eq:noise-variance}, $B$ is a Hilbert--Schmidt operator, with
\begin{equation}\label{eq:B-HS}
    \norm{B}_{L_2(\ell^2,L^2(D))}^2=\sum_n\norm{Bg_n}_{L^2(D)}^2=\sigma^2\sum_n\norm{\zeta_n}_{L^2(D)}^2=\sigma^2.
\end{equation}
By the definition of the stochastic integral with respect to $\mathcal W$ (see \cite[Section~2.5.2]{LR15}), the It\^o isometry and \eqref{eq:B-def}, the series $\sigma\sum_n\zeta_nW^n_t$ of the noise in \eqref{eq:theta-equation} converges in $L^2(\Omega,\mathbb P;L^2(D))$ for every $t\geq0$, and its sum is the stochastic integral of $B$ with respect to $\mathcal W$:
\begin{equation}\label{eq:LR-noise}
    \sigma\sum_n\zeta_nW^n_t=\int_0^tB\,d\mathcal W(s).
\end{equation}

\subsection{Preliminaries on variational solutions}\label{subsec:variational}

We recall an existence theorem for variational solutions and It\^o's formula for the squared norm, stated for the spaces and the noise of this paper. Let $H^1_0(D)$ and $H^{-1}(D)$ be as in \eqref{eq:H1-fourier} and \eqref{eq:Hminus1-fourier}, with norm $\norm{\theta}_{H^1_0(D)}=\norm{\nabla\theta}_{L^2(D)}$ on $H^1_0(D)$, and with the pairing $\langle\cdot,\cdot\rangle_{H^{-1}\times H^1_0}$ between $H^{-1}(D)$ and $H^1_0(D)$. The space $H^1_0(D)$ is a separable Hilbert space, dense in $L^2(D)$, and consequently $H^1_0(D)\subset L^2(D)\subset H^{-1}(D)$ is a Gelfand triple, in which $\langle f,\vartheta\rangle_{H^{-1}\times H^1_0}=(f,\vartheta)$ for $f\in L^2(D)$ and $\vartheta\in H^1_0(D)$.

Let $\mathcal W$ be the cylindrical Wiener process \eqref{eq:cyl-W}, with the filtration $\mathcal F_t$ for which the $W^n$ are Brownian motions, and let $B$ be the operator \eqref{eq:B-def}.

\begin{theorem}[{\cite[Theorem~5.1.3]{LR15}}]\label{thm:LR}
Let $A:H^1_0(D)\to H^{-1}(D)$. Suppose that there is a constant $C_0>0$ such that the following conditions hold for all $v_1,v_2,\vartheta\in H^1_0(D)$:
\begin{enumerate}
    \item[(H1)] (Hemicontinuity) The map $\lambda\mapsto\langle A(v_1+\lambda v_2),\vartheta\rangle_{H^{-1}\times H^1_0}$ is continuous on $\R$.
    \item[(H2$'$)] (Local monotonicity)
    \[
        2\langle A(v_1)-A(v_2),v_1-v_2\rangle_{H^{-1}\times H^1_0}\leq C_0\bigl(1+\norm{v_2}_{L^2(D)}^2\bigr)\norm{v_1-v_2}_{L^2(D)}^2.
    \]
    \item[(H3)] (Coercivity)
    \[
        2\langle A(v_2),v_2\rangle_{H^{-1}\times H^1_0}+\norm{B}_{L_2(\ell^2,L^2(D))}^2\leq C_0\norm{v_2}_{L^2(D)}^2-2\norm{\nabla v_2}_{L^2(D)}^2+C_0.
    \]
    \item[(H4$'$)] (Growth)
    \[
        \norm{A(v_2)}_{H^{-1}(D)}^2\leq C_0\bigl(1+\norm{\nabla v_2}_{L^2(D)}^2\bigr)\bigl(1+\norm{v_2}_{L^2(D)}^2\bigr).
    \]
\end{enumerate}
Then, for every deterministic $\theta_0\in L^2(D)$ and every $t_1>0$, there is a unique $\mathcal F_t$-adapted process $\theta:[0,t_1]\times\Omega\to L^2(D)$ with continuous paths and
\[
    \E\int_0^{t_1}\norm{\nabla\theta(t)}_{L^2(D)}^2\,dt<\infty
\]
such that, almost surely,
\begin{equation}\label{eq:LR-equation}
    \theta(t)=\theta_0+\int_0^tA(\theta(r))\,dr+\int_0^tB\,d\mathcal W(r)\quad\text{in }H^{-1}(D),\qquad\text{for all }t\in[0,t_1].
\end{equation}
This process satisfies
\[
    \E\sup_{t\leq t_1}\norm{\theta(t)}_{L^2(D)}^4<\infty,
\]
and $\theta$ is a Markov process.
\end{theorem}

\begin{remark}\label{rem:LR-noise}
In the notation of \cite{LR15}, Theorem~\ref{thm:LR} is \cite[Theorem~5.1.3]{LR15} for the Gelfand triple $V=H^1_0(D)$, $H=L^2(D)$, $V^*=H^{-1}(D)$, with $U=\ell^2$, $T=t_1$, $W=\mathcal W$, $X_0=\theta_0$, with $A$ and $B$ independent of $t$ and $\omega$, with $B$ independent of the solution, and with
\[
    \alpha=2,\qquad\beta=2,\qquad\theta=2,\qquad p=4,\qquad f\equiv C_0,\qquad\rho(v)=C_0\norm{v}_H^2.
\]
For these values, the conditions (H1), (H2$'$), (H3) and (H4$'$) of \cite[Section~5.1]{LR15} are the conditions of Theorem~\ref{thm:LR}, the two conditions \cite[(5.2)]{LR15} hold because $B$ is constant and $\rho(v)\leq C_0\bigl(1+\norm{v}_V^2\bigr)\bigl(1+\norm{v}_H^2\bigr)$, and a solution in the sense of \cite[Definition~5.1.2]{LR15} is a process $\theta$ as in Theorem~\ref{thm:LR} that satisfies \eqref{eq:LR-equation}.
\end{remark}

We also recall It\^o's formula for the squared norm, stated for the same spaces and for the cylindrical Wiener process $\mathcal W$ in \eqref{eq:cyl-W}.

\begin{theorem}[{\cite[Theorem~4.2.5]{LR15}}]\label{thm:ito-norm}
Let $t_1>0$, let $\xi_0\in L^2(D)$, and let
\[
    h\in L^2\bigl([0,t_1]\times\Omega,dt\otimes\mathbb P;H^{-1}(D)\bigr),\qquad Z\in L^2\bigl([0,t_1]\times\Omega,dt\otimes\mathbb P;L_2(\ell^2,L^2(D))\bigr)
\]
be progressively measurable. Suppose that the continuous $H^{-1}(D)$-valued process
\begin{equation}\label{eq:xi-process}
    \xi(t)\coloneqq\xi_0+\int_0^th(s)\,ds+\int_0^tZ(s)\,d\mathcal W(s),\qquad t\in[0,t_1],
\end{equation}
belongs to $L^2\bigl([0,t_1]\times\Omega,dt\otimes\mathbb P;H^1_0(D)\bigr)$ and satisfies $\E\norm{\xi(t)}_{L^2(D)}^2<\infty$ for almost every $t\in[0,t_1]$. Then $\xi$ is a continuous $L^2(D)$-valued process, and, almost surely, for all $t\in[0,t_1]$,
\begin{equation}\label{eq:ito-norm-formula}
\begin{aligned}
    \norm{\xi(t)}_{L^2(D)}^2&=\norm{\xi_0}_{L^2(D)}^2+\int_0^t\Bigl(2\langle h(s),\xi(s)\rangle_{H^{-1}\times H^1_0}+\norm{Z(s)}_{L_2(\ell^2,L^2(D))}^2\Bigr)\,ds\\
    &\quad+2\int_0^t\bigl(\xi(s),Z(s)\,d\mathcal W(s)\bigr).
\end{aligned}
\end{equation}
\end{theorem}

\begin{remark}\label{rem:ito-norm}
In the notation of \cite[Theorem~4.2.5]{LR15}, Theorem~\ref{thm:ito-norm} is the case $V=H^1_0(D)$, $H=L^2(D)$, $U=\ell^2$, $T=t_1$, $W=\mathcal W$, $\alpha=2$, with $X=\xi$, $X_0=\xi_0$ and $Y=h$.
\end{remark}

\subsection{\texorpdfstring{Proof of Proposition~\ref{prop:wellposed}}{Proof of Proposition 3.2}}\label{subsec:proof-prop-wellposed}

The drift in \eqref{eq:weak-solution} can be viewed as the map $A:H^1_0(D)\to H^{-1}(D)$ given by
\begin{equation}\label{eq:A-drift}
    \langle A(v),\vartheta\rangle_{H^{-1}\times H^1_0}\coloneqq-(\nabla v,\nabla\vartheta)+(v,\vu[v]\cdot\nabla\vartheta)+\Ray(Kv,\vartheta).
\end{equation}
By \eqref{eq:LR-noise}, the equation \eqref{eq:LR-equation} for this map $A$, tested against $\vartheta\in H^1_0(D)$, is \eqref{eq:weak-solution} for $t\in[0,t_1]$.

\paragraph{Existence and proof of (i).}
We wish to show that the map $A$ in \eqref{eq:A-drift} satisfies the conditions (H1), (H2$'$), (H3) and (H4$'$) of Theorem~\ref{thm:LR}, for a sufficiently large constant $C_0$ that depends only on $\Ray$ and $L$.

\emph{Condition (H1).} For $v_1,v_2,\vartheta\in H^1_0(D)$, the map $\lambda\mapsto\langle A(v_1+\lambda v_2),\vartheta\rangle_{H^{-1}\times H^1_0}$ is a polynomial, and hence continuous.

\emph{Condition (H2$'$).} By \eqref{eq:u-bound},
\begin{equation}\label{eq:u-Linf}
    \norm{\vu[v]}_{L^\infty(D)}\leq C_{\vu}\norm{v}_{L^2(D)},
\end{equation}
with a fixed constant $C_{\vu}>0$ depending only on $\Ray$ and $L$. Let $v,\vartheta\in H^1_0(D)$. Since $\divg\vu[v]=0$ and $u_2[v]=0$ at $z=0,1$,
\begin{equation}\label{eq:skew}
    (\vartheta,\vu[v]\cdot\nabla\vartheta)=0.
\end{equation}
Let $v_1,v_2\in H^1_0(D)$ and $w\coloneqq v_1-v_2$. Since $\vu[v]$ is linear in $v$, the identity \eqref{eq:skew} with $v=v_1$ and $\vartheta=w$ implies
\begin{equation}\label{eq:nonlinear-difference}
\begin{aligned}
    (v_1,\vu[v_1]\cdot\nabla w)-(v_2,\vu[v_2]\cdot\nabla w)&=(w,\vu[v_1]\cdot\nabla w)+(v_2,\vu[w]\cdot\nabla w)\\
    &=(v_2,\vu[w]\cdot\nabla w).
\end{aligned}
\end{equation}
By the definition \eqref{eq:A-drift} of $A$ and the identity \eqref{eq:nonlinear-difference},
\[
    2\langle A(v_1)-A(v_2),w\rangle_{H^{-1}\times H^1_0}=-2\norm{\nabla w}_{L^2(D)}^2+2(v_2,\vu[w]\cdot\nabla w)+2\Ray(Kw,w).
\]
Therefore, the Cauchy--Schwarz inequality, the estimate \eqref{eq:u-Linf} for $\vu[w]$, the estimate \eqref{eq:K-eig-bound} on the eigenvalues of $K$, and Young's inequality $2ab\leq2a^2+\frac12b^2$ imply
\begin{equation}\label{eq:drift-monotone}
\begin{aligned}
    2\langle A(v_1)-A(v_2),w\rangle_{H^{-1}\times H^1_0}&\leq-2\norm{\nabla w}_{L^2(D)}^2+2C_{\vu}\norm{v_2}_{L^2(D)}\norm{w}_{L^2(D)}\norm{\nabla w}_{L^2(D)}\\
    &\qquad+\frac{\Ray}{2\pi^2}\norm{w}_{L^2(D)}^2\\
    &\leq\Bigl(\frac{\Ray}{2\pi^2}+\frac{C_{\vu}^2}2\norm{v_2}_{L^2(D)}^2\Bigr)\norm{w}_{L^2(D)}^2.
\end{aligned}
\end{equation}
The estimate \eqref{eq:drift-monotone} implies (H2$'$).

\emph{Condition (H3).} For $v\in H^1_0(D)$, the definition \eqref{eq:A-drift} of $A$, the identity \eqref{eq:skew} with $\vartheta=v$, and the estimate \eqref{eq:u2-theta-bound} imply
\begin{equation}\label{eq:drift-coercive}
    2\langle A(v),v\rangle_{H^{-1}\times H^1_0}=-2\norm{\nabla v}_{L^2(D)}^2+2(u_2[v],v)\leq\frac{\Ray}{2\pi^2}\norm{v}_{L^2(D)}^2-2\norm{\nabla v}_{L^2(D)}^2.
\end{equation}
Since $\norm{B}_{L_2(\ell^2,L^2(D))}^2=\sigma^2$ by \eqref{eq:B-HS} and $\sigma\leq1$, the estimate \eqref{eq:drift-coercive} implies (H3).

\emph{Condition (H4$'$).} The definition \eqref{eq:A-drift} of $A$, the Cauchy--Schwarz inequality, the estimates \eqref{eq:u-Linf} and \eqref{eq:K-eig-bound}, and Poincar\'e's inequality $\norm{\vartheta}_{L^2(D)}\leq\frac1\pi\norm{\nabla\vartheta}_{L^2(D)}$ for $\vartheta\in H^1_0(D)$ imply
\begin{equation}\label{eq:drift-growth}
    \norm{A(v)}_{H^{-1}(D)}\leq\norm{\nabla v}_{L^2(D)}+C_{\vu}\norm{v}_{L^2(D)}^2+\frac{\Ray}{4\pi^3}\norm{v}_{L^2(D)}.
\end{equation}
Squaring \eqref{eq:drift-growth} and using $\norm{v}_{L^2(D)}^4\leq\pi^{-2}\norm{\nabla v}_{L^2(D)}^2\norm{v}_{L^2(D)}^2$, we obtain (H4$'$).

Consequently, Theorem~\ref{thm:LR} establishes, for every $t_1>0$, a solution on $[0,t_1]$ of \eqref{eq:LR-equation}, with
\begin{equation}\label{eq:LR-moments}
    \E\sup_{r\leq t_1}\norm{\theta(r)}_{L^2(D)}^4+\E\int_0^{t_1}\norm{\nabla\theta(r)}_{L^2(D)}^2\,dr<\infty.
\end{equation}
By the uniqueness in Theorem~\ref{thm:LR}, the solutions on $[0,t_1]$ for different $t_1>0$ agree on their common time interval, and therefore define a process $\theta$ on $[0,\infty)$. This process satisfies \eqref{eq:weak-solution}, by \eqref{eq:LR-noise}, and the moment bound \eqref{eq:solution-moments}, by \eqref{eq:LR-moments}, and therefore $\theta$ is a solution in the sense of Definition~\ref{def:solution}. By the uniqueness in Theorem~\ref{thm:LR}, $\theta$ is the unique solution. Moreover, \eqref{eq:LR-moments} implies the estimate in (i). Finally, Theorem~\ref{thm:LR} states that $\theta$ is a Markov process, and this Markov process is time-homogeneous (see \cite[Proposition~4.3.5]{LR15}), as (i) states.

\paragraph{Proof of (ii).}
We follow the proof of \cite[Lemma~5.2]{Sch17}. For $t\geq0$, $y\in L^2(D)$, Borel sets $E\subset L^2(D)$ and bounded Borel functions $\varphi:L^2(D)\to\R$, denote
\[
    P_t(y,E)\coloneqq\mathbb P_y(\theta(t)\in E),\qquad P_t\varphi(y)\coloneqq\E_y\varphi(\theta(t))=\int_{L^2(D)}\varphi(y')\,P_t(y,dy').
\]
We wish to show the joint measurability of the transition kernel $(t,y)\mapsto P_t(y,E)$.

We start with an estimate on the difference of two solutions. Let $\theta^{(1)},\theta^{(2)}$ be solutions in the sense of Definition~\ref{def:solution}, for the same Brownian motions $W^n$, with any initial values $\theta^{(1)}(0),\theta^{(2)}(0)\in L^2(D)$, and let $w\coloneqq\theta^{(1)}-\theta^{(2)}$. Subtracting the identities \eqref{eq:weak-solution} for $\theta^{(1)}$ and $\theta^{(2)}$, the noise terms $\sigma\sum_n(\zeta_n,\vartheta)W^n_t$ cancel, and, almost surely, for all $t$,
\[
    w(t)=w(0)+\int_0^t\bigl(A(\theta^{(1)}(r))-A(\theta^{(2)}(r))\bigr)\,dr\qquad\text{in }H^{-1}(D).
\]
By \eqref{eq:drift-growth} and the moment bound \eqref{eq:LR-moments}, $A(\theta^{(1)})$ and $A(\theta^{(2)})$ belong to $L^2([0,t]\times\Omega,dt\otimes\mathbb P;H^{-1}(D))$, $w$ belongs to $L^2([0,t]\times\Omega,dt\otimes\mathbb P;H^1_0(D))$ and $\E\norm{w(r)}_{L^2(D)}^2<\infty$ for every $r\leq t$. Therefore, Theorem~\ref{thm:ito-norm}, with $\xi=w$, $\xi_0=w(0)$, $h=A(\theta^{(1)})-A(\theta^{(2)})$ and $Z=0$, and \eqref{eq:drift-monotone} yield
\[
    \norm{w(t)}_{L^2(D)}^2\leq\norm{w(0)}_{L^2(D)}^2+\int_0^t\Bigl(\frac{\Ray}{2\pi^2}+\frac{C_{\vu}^2}2\norm{\theta^{(2)}(r)}_{L^2(D)}^2\Bigr)\norm{w(r)}_{L^2(D)}^2\,dr,
\]
and from Gr\"onwall's inequality we have
\begin{equation}\label{eq:theta-stability}
    \norm{w(t)}_{L^2(D)}^2\leq\norm{w(0)}_{L^2(D)}^2\exp\Bigl(\frac{\Ray t}{2\pi^2}+\frac{C_{\vu}^2}2\int_0^t\norm{\theta^{(2)}(r)}_{L^2(D)}^2\,dr\Bigr).
\end{equation}

Notice that $(t,y)\mapsto P_t\varphi(y)$ is continuous for every bounded continuous $\varphi:L^2(D)\to\R$. This follows from
\[
    \abs{P_{t_m}\varphi(y_m)-P_t\varphi(y)}\leq\abs{P_{t_m}\varphi(y_m)-P_{t_m}\varphi(y)}+\abs{P_{t_m}\varphi(y)-P_t\varphi(y)},
\]
for $y_m\to y$ in $L^2(D)$ and $t_m\to t$, the estimate \eqref{eq:theta-stability} and the continuity of the paths of $\theta$. Indeed, for each $m\geq1$, let $\theta^m$ be the solution with $\theta^m(0)=y_m$, and let $\theta$ be the solution with $\theta(0)=y$, all for the same Brownian motions $W^n$. By \eqref{eq:theta-stability} with $\theta^{(1)}=\theta^m$ and $\theta^{(2)}=\theta$, we have $\theta^m(t_m)-\theta(t_m)\to0$ in $L^2(D)$, almost surely, and $\theta(t_m)\to\theta(t)$ in $L^2(D)$, because the paths of $\theta$ are continuous. Therefore, both terms on the right side tend to $0$ as $m\to\infty$, by the dominated convergence theorem.

For every nonempty closed set $E\subset L^2(D)$, the functions defined by
\[
    \varphi_n(y)\coloneqq\max\{0,\,1-n\,d(y,E)\},\qquad d(y,E)\coloneqq\inf_{e\in E}\norm{y-e}_{L^2(D)},
\]
are bounded and continuous, and $\varphi_n\downarrow\mathbf 1_E$. By monotone convergence, $P_t(y,E)=\inf_{n\geq1}P_t\varphi_n(y)$, which proves the joint measurability in $(t,y)$ for all closed sets.

Since
\[
    \mathcal D\coloneqq\bigl\{E\subset L^2(D)\text{ Borel}:\ (t,y)\mapsto P_t(y,E)\text{ is measurable}\bigr\}
\]
is a Dynkin system containing the closed sets, which are stable under finite intersections and generate the Borel $\sigma$-algebra of $L^2(D)$, $\mathcal D$ contains every Borel subset of $L^2(D)$.

\paragraph{Proof of (iii).}
Let $\theta$ be the solution with $\theta(0)=\theta_0$, and let $t_1>0$. By \eqref{eq:LR-equation}, $\theta$ is the process \eqref{eq:xi-process} with $\xi_0=\theta_0$, $h=A(\theta)$ and $Z=B$, the operator \eqref{eq:B-def}. By \eqref{eq:drift-growth}, (i) and \eqref{eq:solution-moments}, the hypotheses of Theorem~\ref{thm:ito-norm} hold, and \eqref{eq:ito-norm-formula}, with \eqref{eq:B-HS}, \eqref{eq:LR-noise} and the identity in \eqref{eq:drift-coercive}, yields \eqref{eq:energy-ito}. By (i) and $\sum_n(\theta,\zeta_n)^2\leq\norm{\theta}_{L^2(D)}^2$, the stochastic integral is a martingale.

\section*{Statement on the use of artificial intelligence}

During the preparation of this work the authors used the large language model Claude Fable 5.2. The model pointed out that Appendix~\ref{sec:wellposed} has to write the noise in \eqref{eq:theta-equation} as a stochastic integral with respect to a cylindrical Wiener process, which Section~\ref{subsec:noise-cylindrical} does. The model also pointed out that the constants in the proof of Proposition~\ref{prop:construct-Y} have to be chosen in a careful order to avoid a contradiction, which that proof does. The authors take full responsibility for the content of this article.

\end{document}